\batchmode
\documentclass[11pt]{article}

\usepackage{epsfig}
\usepackage{graphicx}
\usepackage{color}
\usepackage{mathtools}

\newtheorem{theorem}{Theorem}
\newtheorem{lemma}{Lemma}
\newtheorem{corollary}{Corollary}

\newtheorem{definition}{Definition}

\newcommand{\be}{\begin{equation}}
\newcommand{\ee}{\end{equation}}
\newcommand{\bea}{\begin{eqnarray}}
\newcommand{\eea}{\end{eqnarray}}
\newcommand{\beas}{\begin{eqnarray*}}
\newcommand{\eeas}{\end{eqnarray*}}
\newcommand{\ba}{\begin{array}}
\newcommand{\ea}{\end{array}}

\DeclarePairedDelimiter{\floor}{\lfloor}{\rfloor}
\DeclarePairedDelimiter{\ceil}{\lceil}{\rceil}

\definecolor{armygreen}{rgb}{0.29, 0.33, 0.13}

\newcommand{\real}{\mbox{$\mathbb{R}$}}

\newcommand{\bfb}{\ensuremath{\mathbf{b}}}
\newcommand{\bfc}{\ensuremath{\mathbf{c}}}

\def\XXint#1#2#3{{\setbox0=\hbox{$#1{#2#3}{\int}$}
     \vcenter{\hbox{$#2#3$}}\kern-.5\wd0}}

\newcommand{\mcI}{\ensuremath{\mathcal{I}}}

\newcommand{\mcL}{\ensuremath{\mathcal{L}}}

\newcommand{\mcS}{\ensuremath{\mathcal{S}}}

\newcommand{\mrI}{\ensuremath{\mathrm{I}}}
\newcommand{\mrJ}{\ensuremath{\mathrm{J}}}

\newcommand{\wtilde}{\ensuremath{\widetilde}}
\newcommand{\what}{\ensuremath{\widehat}}

\def\qed{\hbox{\vrule width 6pt height 6pt depth 0pt}}

\usepackage{amsmath}
\usepackage{amssymb}

\title{Optimal Petrov-Galerkin spectral approximation method for the fractional diffusion, 
advection, reaction equation on a bounded interval} 
\author{
	Xiangcheng Zheng \thanks{Department of Mathematics, University of South Carolina, Columbia,
		South Carolina 29208, USA. email: {\tt xz3@math.sc.edu \& hwang@math.sc.edu}.} 
	\and
	V.J.~Ervin\thanks{School of Mathematical and Statistical Sciences,
	  Clemson University, Clemson, South Carolina 29634-0975, USA.
	  email: {\tt vjervin@clemson.edu}. }
	\and 
	 Hong Wang $^*$ }

\date{\today}

\begin{document}
\maketitle

\begin{abstract}
In this paper we investigate the numerical approximation of the fractional diffusion, advection, reaction equation
on a bounded interval.
Recently the explicit form of the solution to this equation was obtained. Using the explicit form of
the boundary behavior of the solution and Jacobi polynomials, a Petrov-Galerkin approximation
scheme is proposed and analyzed. Numerical experiments are presented which support the
theoretical results, and demonstrate the accuracy and optimal convergence of the approximation method.
\end{abstract}



\setcounter{equation}{0}
\setcounter{figure}{0}
\setcounter{table}{0}
\setcounter{theorem}{0}
\setcounter{lemma}{0}
\setcounter{corollary}{0}
\setcounter{definition}{0}
%
\section{Introduction}
 \label{sec_intro}
Of interest in this paper is the approximation of the solution to the fractional diffusion, 
advection, reaction equation
\begin{align}
 \mcL_{r}^{\alpha}u(x) \ + \ b(x) D u(x) \ + \ c(x) u(x)   &= \ f(x) \, , \ \  x \in \mrI \, , \
   \label{DefProb2}  \\
   \mbox{subject to } u(0) \, = \, u(1) &= \,  0 \, ,
\label{DefBC2}  \\
\mbox{where } \ \mcL_{r}^{\alpha}u(x) \ := \  
   - D \big( r D^{-(2 - \alpha)} &+ \ (1 - r) D^{-(2 - \alpha)*}  \big) D u(x)  ,   \label{defmcL}
\end{align}  
and $\mrI \, := \, (0 , 1)$, $1 < \alpha < 2$, $0 \le r \le 1$, $c(x) \, - \, \frac{1}{2}D b(x) \ge 0$, 
$D$ denotes the usual derivative
operator, $D^{\alpha}$ the $\alpha$-order left fractional derivative operator, and
$D^{\alpha *}$ the $\alpha$-order right fractional derivative operator, defined by:
\begin{align}
D^{\alpha} u(x) &:= \ D \, \mbox{}_{0}D_{x}^{- (2 - \alpha)} \, D u(x) \ 
= \ D \, \frac{1}{\Gamma(2 - \alpha)}   \int_{0}^{x} \frac{1}{(x - s)^{\alpha - 1}} \, D u(s) \, ds \, ,   \label{defDalpha}  \\
D^{\alpha *} u(x) &:= \ D \, \mbox{}_{x}D_{1}^{- (2 - \alpha)} \, D u(x) \ 
= \ D \,  \frac{1}{\Gamma(2 - \alpha)} \int_{x}^{1} \frac{1}{(s - x)^{\alpha - 1}} \, D u(s) \, ds \, .   \label{defDalpha*} 
\end{align}

In recent years fractional differential equations have received increased attention as they have been used in modeling
a number of physical phenomena such as  contaminant transport in ground water flow \cite{ben001},
viscoelasticity \cite{mai971}, image processing \cite{bua101, gat151},
turbulent flow \cite{mai971, shl871}, and chaotic dynamics \cite{zas931}.

The are two important properties that distinguish a fractional order differential equations from its integer order
counterpart. Firstly, as can be noted from \eqref{defmcL}, fractional differential equations are nonlocal in nature. Secondly, 
the solution of fractional differential equations (typically) have a lack of regularity at the boundary of the domain.
Finite difference methods \cite{cui091, liu041, mee041, tad071, wan121},
finite element methods \cite{erv061, jin151, liu111, wan131}, 
discontinuous Galerkin methods \cite{xu141}, and mixed methods \cite{che161, li171}, have all been developed for 
fractional differential equations. These methods typically exhibit slow convergence due to the lack of regularity
of the solution at the boundary. In \cite{jin161, jin162} an enriched subspace was given for one sided fractional 
differential equations, where the boundary behavior of the solution was included in the finite element trial space.
Mao and Shen in \cite{mao182} extended the work of Gui and Babu\v{s}ka in \cite{gui861} to establish that
for an assumed boundary behavior of the solution a geometrically spaced mesh with increasing polynomial
degree trial function on the subintervals resulted in an exponential rate of converge for the approximation.
For a special class of self-adjoint fractional differential equations a spectral approximation scheme was 
presented in \cite{zay131} using a special class of functions, polyfractonomials. Spectral 
methods, exploiting a special property satisfied by fractional diffusion operator applied to Jacobi polynomials
(see \eqref{propmcL}) has been particularly effective for the approximation of the solution to fractional 
diffusion equations \cite{che162, erv162, li121, mao161, mao162, mao181, zhe191, zhe201}. 

Two recent papers have established the explicit form of solutions to fractional diffusion, advection, reaction equations on
a bounded domain in $\mathbb{R}^{1}$. In \cite{hao201}, Hao and Zhang studied the case for $r = 1/2$, for which 
$\mcL_{r}^{\alpha}$ is a symmetric operator. Their work was extended in \cite{erv191} to the general case of $0 \le r \le 1$.
The solution was shown to have the form $u(x) \ = \ (1 - x)^{\alpha - \beta} x^{\beta} \phi(x)$, where 
$\phi$ is contained in the weighted Sobolev space $H^{\alpha \, + \, \wtilde{s}}_{(\alpha - \beta \, , \, \beta)}(\mrI)$ (defined
in Section \ref{sec_not}), where $\beta$ and $\wtilde{s}$ are explicit functions of $\alpha, \, r$, and the
regularity of the right hand side function, $f$ (see Theorems \ref{thmreg11} and \ref{thmreg13} below). 
Of particular note is that for the fractional diffusion,
reaction problem, and the fractional diffusion, advection, reaction problem, the regularity of the solution $u$ is
bounded, regardless of the regularity of $f$. This boundedness in the regularity of $u$ is not the case for the 
fractional diffusion, advection, reaction equation on $\mathbb{R}$, as was recently established by Ginting and Li
in \cite{gin191}.

The numerical approximation scheme presented below is accurate as, using \cite{erv191}, the precise
boundary behavior of the solution is incorporated into the approximate solution. Additionally, using the special property
of the fractional diffusion operator applied to Jacobi polynomials (see \eqref{propmcL})
\[
 \mcL^{\alpha}_{r} \omega(x) \, \what{G}_{k}^{(\alpha - \beta \, , \, \beta)}(x)
   \ = \ \lambda_{k} \, \what{G}_{k}^{(\beta \, , \, \alpha - \beta)}(x) \, , 
\]
and that $\{   \what{G}_{k}^{(\alpha - \beta \, , \, \beta)} \}_{k = 0}^{\infty}$ is a  basis for 
$H^{r}_{(\alpha - \beta \, , \, \beta)}(\mrI)$, the approximation scheme using Jacobi polynomial is efficient in that
if the solution is $C^{\infty}(\mrI)$ (very rarely the case) the approximation converges exponentially. If the solution 
has bounded regularity (typically the case) the approximation converges optimally at an algebraic rate of convergence.

This paper is organized as follows. In the following section definitions, notation, and several known results
are summarized. Section \ref{sec_form} contains the Petrov-Galerkin weak formulation for
\eqref{DefProb2},\eqref{DefBC2}, and establishes the existence and uniqueness of its solution.
The analysis follows the work of Jin, Lazarov and Zhou in \cite{jin162}, wherein the 
lower order terms are handled using the Petree-Tartar Lemma.
The approximation scheme is given in Section \ref{sec_appx}, and associated error estimates derived.
Numerical experiments are presented in Section \ref{sec_num}.

 \setcounter{equation}{0}
\setcounter{figure}{0}
\setcounter{table}{0}
\setcounter{theorem}{0}
\setcounter{lemma}{0}
\setcounter{corollary}{0}
\setcounter{definition}{0}
%
\section{Notation and Properties}
\label{sec_not}
Jacobi polynomials have an important connection with fractional order diffusion equations
\cite{aco181, erv162, mao181, mao161}. We briefly review their definition and some of their important
properties \cite{abr641, sze751}. 

\textbf{Usual Jacobi Polynomials, $P_{n}^{(a , b)}(t)$, on $(-1 \, , \, 1)$}.   \\    
\textbf{Definition}: $ P_{n}^{(a , b)}(t) \ := \ 
\sum_{m = 0}^{n} \, p_{n , m} \, (t - 1)^{(n - m)} (t + 1)^{m}$, where
\begin{equation}
       p_{n , m} \ := \ \frac{1}{2^{n}} \, \left( \begin{array}{c}
                                                              n + a \\
                                                              m  \end{array} \right) \,
                                                    \left( \begin{array}{c}
                                                              n + b \\
                                                              n - m  \end{array} \right) \, .
\label{spm21}
\end{equation}
\underline{Orthogonality}:    
\begin{align}
 & \int_{-1}^{1} (1 - t)^{a} (1 + t)^{b} \, P_{j}^{(a , b)}(t) \, P_{k}^{(a , b)}(t)  \, dt 
 \ = \
   \left\{ \begin{array}{ll} 
   0 , & k \ne j  \\
   |\| P_{j}^{(a , b)} |\|^{2}
   \, , & k = j  
    \end{array} \right.  \, ,  \nonumber \\
& \quad \quad \mbox{where } \  \ |\| P_{j}^{(a , b)} |\| \ = \
 \left( \frac{2^{(a + b + 1)}}{(2j \, + \, a \, + \, b \, + 1)} 
   \frac{\Gamma(j + a + 1) \, \Gamma(j + b + 1)}{\Gamma(j + 1) \, \Gamma(j + a + b + 1)}
   \right)^{1/2} \, .
  \label{spm22}
\end{align}                                                    

In order to transform the domain of the family of Jacobi polynomials to $[0 , 1]$, let $t \rightarrow 2x - 1$ and 
introduce $G_{n}^{(a , b)}(x) \, = \, P_{n}^{(a , b)}( t(x) )$. From \eqref{spm22},
\begin{align}
 \int_{-1}^{1} (1 - t)^{a} (1 + t)^{b} \, P_{j}^{(a , b)}(t) \, P_{k}^{(a , b)}(t)  \, dt 
 &= \
 \int_{0}^{1} 2^{a} \, (1 - x)^{a} \, 2^{b} \, x^{b} \, P_{j}^{(a , b)}(2x - 1) \, 
 P_{k}^{(a , b)}(2x - 1)  \, 2 \,  dx
 \nonumber \\
  &= \
2^{a + b + 1} \int_{0}^{1}   (1 - x)^{a}  \, x^{b} \, G_{j}^{(a , b)}(x) \, G_{k}^{(a , b)}(x)  \,  dx
\nonumber \\
&= \
   \left\{ \begin{array}{ll} 
   0 , & k \ne j \, , \\
  2^{a + b + 1} \, |\| G_{j}^{(a , b)} |\|^{2}
   \, , & k = j  
   \, . \end{array} \right.    \nonumber \\
 \quad \quad \mbox{where } \  \ |\| G_{j}^{(a , b)} |\| &= \
 \left( \frac{1}{(2j \, + \, a \, + \, b \, + 1)} 
   \frac{\Gamma(j + a + 1) \, \Gamma(j + b + 1)}{\Gamma(j + 1) \, \Gamma(j + a + b + 1)}
   \right)^{1/2} \, .  \label{spm22g} 
\end{align}                                                    

 From \cite[equation (2.19)]{mao161} we have that
\begin{equation}
   \frac{d^{k}}{dt^{k}} P_{n}^{(a , b)}(t) \ = \ 
   \frac{\Gamma(n + k + a + b + 1)}{2^{k} \, \Gamma(n + a + b + 1)} P_{n - k}^{(a + k \, , \, b + k)}(t) \, .
   \label{derP}
\end{equation}   
Hence,
\begin{align}
\frac{d^{k}}{dx^{k}} G_{n}^{(a , b)}(x) 
  &= \ \frac{\Gamma(n + k + a + b + 1)}{  \Gamma(n + a + b + 1)} 
  G_{n - k}^{(a + k \, , \, b + k)}(x)  \, .  \label{eqC4}
\end{align}

Note that, from Stirling's formula, we have that
\begin{equation}
\lim_{n \rightarrow \infty} \, \frac{\Gamma(n + \sigma)}{\Gamma(n) \, n^{\sigma}}
\ = \ 1 \, , \mbox{ for } \sigma \in \mathbb{R}.  
 \label{eqStrf}
\end{equation} 
 
For compactness of notation, let
\begin{equation}
 \omega^{(a , b)} \, = \, \omega^{(a , b)}(x) \, := \, (1 - x)^{a} \, x^{b} \, .
 \label{defrho}
\end{equation} 

 We let $\mathbb{N}_{0}  := \mathbb{N} \cup {0}$ and
 use $y_{n} \sim n^{p}$ to denote that there exists constants $c$ and $C > 0$ such that, as 
 $n \rightarrow \infty$,  
 $c \, n^{p} \le | y_{n} | \le C \, n^{p}$. Additionally, we use $a \, \lesssim \, b$ to denote that there exists a constant $C$ such that
 $a \, \le \, C \,  b$. 
 
 For $t \in \mathbb{R}$, $\floor{t}$ is used to denote the largest integer that is less than or equal to $t$, and
 $\ceil{t}$ is used to denote the smallest integer that is greater than or equal to $t$.

\textbf{Function space $L_{\sigma}^{2}(\mrI)$}. \\
For $\sigma(x) > 0, \ x \in (0 , 1)$, let 
\begin{equation}
L_{\sigma}^{2}(\mrI) \, := \, \{ f(x) \, : \, \int_{0}^{1} \sigma(x) \, f(x)^{2} \, dx \ < \ \infty \} \, .
\label{defLw}
\end{equation}
Associated with $L_{\sigma}^{2}(0 , 1)$ is the inner product, $( \cdot , \cdot )_{\sigma}$, and
norm, $\| \cdot \|_{\sigma}$, defined by
\[
( f \,  , \, g )_{\sigma} \, := \, \int_{0}^{1} \sigma(x) \, f(x) \, g(x) \, dx \, , \quad \mbox{and} \quad
 \| f \|_{\sigma} \, := \, \left( \langle f \,  , \, f \rangle_{\sigma} \right)^{1/2} \, .
\]

The set of orthogonal polynomials $\{ G_{j}^{(a , b)} \}_{j = 0}^{\infty}$ form an orthogonal basis
for $L^{2}_{\omega^{(a , b)}}(\mrI)$, and for $\widehat{G}_{j}^{(a , b)} \, := \, G_{j}^{(a , b)} / |\| G_{j}^{(a , b)} |\|$, 
$\{ \what{G}_{j}^{(a , b)} \}_{j = 0}^{\infty}$ form an orthonormal basis
for $L^{2}_{\omega^{(a , b)}}(\mrI)$.

Without a subscript, $( \cdot , \cdot )$ denotes the usual $L^{2}(\mrI)$ inner product.

\textbf{Function space $H^{s}_{(a , b)}(\mrI)$}. \\
The weighted Sobolev spaces $H^{s}_{(a , b)}(\mrI)$ differ
from the usual $H^{s}(\mrI)$ spaces in that the associated norms apply a polynomial weight at each
endpont of $\mrI$, namely, $x^{b}$ and $(1 - x)^{a}$. These weights increase with the order of the 
derivative. We give two equivalent definitions for the $H^{s}_{(a , b)}(\mrI)$ spaces. In the first
definition the spaces $H^{s}_{(a , b)}(\mrI)$, for $0 < s \not \in \mathbb{N}$, are defined by the
$K$- method of interpolation. The second definition is based on the decay rate of the Jacobi 
coefficients of a function expanded in terms of the Jacobi polynomials $\what{G}_{j}^{(a , b)}(x)$.
Both definitions are useful, and used in the analysis below.
The equivalence of the spaces is discussed in \cite{erv191}.

\textbf{Definition}: Using Interpolation \\
Following Babu\v{s}ka and Guo \cite{bab011}, and Guo and Wang \cite{guo041}, we introduce the
weighted Sobolev spaces $H^s_{\omega^{(a,b)} }(\mrI)$.
\begin{definition} \label{Defspace0}
Let $s, a, b \in \real$, $s \ge 0$, $a, b > -1$. Then
\begin{equation} 
H^s_{\omega^{(a,b)} }(\mrI) := \bigg \{ v \, : \| v \|_{s , \omega^{(a,b)}}^2 := 
 \sum_{j=0}^s \big \| D^j v \big \|_{\omega^{(a+j,b+j)}}^2 < \infty \bigg \} .
\label{defHw}
\end{equation} 
Definition (\ref{defHw}) is extended to $s \in \mathbb{R}^{+}$ using the $K$- method of interpolation. For $s < 0$
the spaces are defined by (weighted) $L^{2}$ duality.
\end{definition}

\textbf{Definition}: Using the decay rate of Jacobi coefficients \\
Next we define function spaces in terms of 
the decay property of the Jacobi coefficients of their member functions. 

Given $v$, let
\be
  v_{j} \ = \ \int_{0}^{1}  \omega^{(a , b)}(x) \, v(x) \,  \what{G}_{j}^{(a , b)}(x) \, dx \, .
 \label{defvj}
\ee
Note that for $v \in L^{2}_{\omega^{(a , b)}}(\mrI)$,
\be
   v(x) \ = \ \sum_{j = 0}^{\infty} v_{j} \, \what{G}_{j}^{(a , b)}(x) \, .
 \label{vexp1}
\ee

\begin{definition} \label{Defspace}
Let $s, a, b \in \real$,  $a, b > -1$, 
$L^{2}_{(a , b)}(\mrI) \, := \, L^{2}_{\omega^{(a , b)}}(\mrI)$, and
$v_{j}$ be given by \eqref{defvj}. Then, define 
\be
H^{s}_{(a , b)}(\mrI) \, := \, \{v \, : \,
 \sum_{j = 0}^{\infty} (1 + j^{2})^{s} \, v_{j}^{2} \, < \, \infty \}
\label{defHr}
\ee
as the $(a , b)$-weighted Sobolev space of order $s$.
\end{definition}

\begin{theorem}  \cite[Theorem 4.1]{erv191} \label{thmeq2}
The spaces $H^{s}_{(a , b)}(\mrI)$ and $H^{s}_{\omega^{(a , b)}}(\mrI)$
coincide, and their corresponding norms are equivalent. 
\end{theorem}

With the structure of the $H^{s}_{(a , b)}(\mrI)$ spaces, and properties \eqref{eqC4} and \eqref{spm22g},
it is straight forward to show that $D$ is a bounded
mapping from $H^{s}_{(a , b)}(\mrI)$ onto
$H^{s - 1}_{(a + 1 \, , \, b + 1)}(\mrI)$.

\begin{lemma}  \cite[Lemma 4.5]{erv191}\label{lmamapD}
For $s, a, b \in \real$,  $a, b > -1$, the differential operator $D$ is a bounded mapping from $H^{s}_{(a , b)}(\mrI)$ onto
$H^{s - 1}_{(a + 1 \, , \, b + 1)}(\mrI)$.
\end{lemma}

For convenience, from  hereon we use $H^{s}_{(a , b)}(\mrI)$ 
to represent the spaces
$H^s_{\omega^{(a,b)} }(\mrI)$ and $H^{s}_{(a , b)}(\mrI)$.

\textbf{Definition}: \textbf{Condition A} \\
The parameters $a$, $b$, and $r$ and constant $c_{*}^{*}$ satisfy:
$1 < \alpha < 2$, $\alpha - 1 \, \le  \, \beta \, , \, \alpha - \beta \, \le  \, 1$, $0 \le r \le 1$ 
\begin{equation}
  c_{*}^{*} \ = \ \frac{\sin(\pi \alpha)}{\sin(\pi (\alpha - \beta)) \, + \, \sin(\pi \beta)} \, ,  \label{defcss}
\end{equation}  
where $\beta$ is determined by
\begin{equation}
  r \ = \ \frac{\sin( \pi \, \beta)}{\sin( \pi ( \alpha - \beta)) \, + \,  \sin( \pi \, \beta)} \, . \label{propker0} 
\end{equation}

For compactness of notation, for $\alpha$ and $r$ defined in \eqref{DefProb2} and $\beta$ defined in 
\eqref{propker0} we introduce
\begin{equation}
 \omega(x) \, := \, \omega^{(\alpha - \beta , \beta)}(x) \, = \, (1 - x)^{\alpha - \beta} \, x^{\beta} \, ,  \ 
 \mbox{ and } \ \omega^{*}(x) \, := \, \omega^{(\beta  , \alpha - \beta)}(x) \, = \, (1 - x)^{\beta} \, x^{\alpha - \beta} \, .  \label{defomega}
\end{equation} 

Additionally, we use $\langle \cdot , \cdot \rangle_{\omega}$ to denote the weighted $L^{2}$ duality
pairing between functions if $H^{-s}_{(\alpha - \beta \, , \, \beta)}(\mrI)$ and 
$H^{s}_{(\alpha - \beta \, , \, \beta)}(\mrI)$.

From \cite{erv162, jia181},
\be
   \mcL^{\alpha}_{r} \omega(x) \, \what{G}_{k}^{(\alpha - \beta \, , \, \beta)}(x)
   \ = \ \lambda_{k} \, \what{G}_{k}^{(\beta \, , \, \alpha - \beta)}(x) \, , \ \ 
   \mbox{ where } \ \ 
   \lambda_{k} \ = \ - c_{*}^{*} \frac{\Gamma(k + 1 + \alpha)}{\Gamma(k + 1)} \, , \ k = 0, 1, 2, \ldots,
\label{propmcL}
\ee
and $c_{*}^{*} $ given by \eqref{defcss}. Also, using \eqref{eqStrf}, $\lambda_{k} \sim k^{\alpha}$.

Let $\mcS_{N}$ denote the space of polynomials of degree less than or equal to $N$. 
We define the weighted $L^2$ orthogonal projection
$P_{N} : \, L^{2}_{\omega}(\mrI) \rightarrow \mcS_{N}$ by the condition
\begin{equation}\label{Proj}
\big ( v \, - \,  P_{N}v \ , \ \phi_N \big)_{\omega} \ = \ 0 \, , \ \ \forall \phi_N \in \mcS_{N}.
\end{equation}
Note that $P_{N}v \ = \ \sum_{j = 0}^{N} v_{j} \, \what{G}_{j}^{(a , b)}(x)$, 
where $v_{j}  \ = \ \int_{0}^{1}  \omega(x) \, v(x) \,  \what{G}_{j}^{(a , b)}(x) \, dx$.
\begin{lemma}\label{lem:Approx} \cite[Theorem 2.1]{guo041}
For $\mu \in \mathbb{N}_{0}$ and $v \in H^{t}_{\omega}(\mrI)$, with $0 \le \mu \le t$, there exists a
constant $C$, independent of $N, \, \alpha$ and $\beta$ such that
\begin{equation}\label{Approx}
\big \| v  -  P_{N} v \|_{H^{\mu}_{\omega}(\mrI)} \ \le \ C \, 
N^{\mu - t} \, \| v \|_{H^{t}_{\omega}(\mrI)}.
\end{equation}
\end{lemma}  

\textbf{Remark}: In \cite{guo041} \eqref{Approx} is stated for $t \in \mathbb{N}_{0}$. The result extends to 
$t \in \mathbb{R}^{+}$ using interpolation.


The regularity of the solution to \eqref{DefProb2}
can be influenced by the regularity of the coefficients $b(x)$ and $c(x)$. The following lemma enables
us to insulate the influence of these terms.

Introduce the space $W^{k , \infty}_{w}(\mrI)$ and its associated norm, defined for $k \in \mathbb{N}_{0}$, as
\begin{align}
W^{k , \infty}_{w}(\mrI) &:= \ \left\{ f \, : \ (1 - x)^{j/2} x^{j/2} D^{j}f(x) \in L^{\infty}(\mrI) , \ j = 0, 1, \ldots, k \right\} , 
\label{defWinftw}  \\
\| f \|_{W^{k , \infty}_{w}} &:= \ \max_{0 \le j \le k} \| (1 - x)^{j/2} x^{j/2} D^{j}f(x) \|_{L^{\infty}(\mrI)} \, .
\label{defWinftwnorm}
\end{align}
The subscript $w$ denotes the fact that $W^{k , \infty}_{w}(\mrI)$ is a weaker space than $W^{k , \infty}(\mrI)$ in that
the derivative of functions in $W^{k , \infty}_{w}(\mrI)$ may be unbounded at the endpoints of the interval. 
\begin{lemma} \label{lmaprodsp}  \cite[Lemma 7.1]{erv191}
Let $s \ge 0$, $\alpha, \, \beta > -1$, $k \ge s$, and $f \in W^{k , \infty}_{w}(\mcI)$. For
\begin{align}
(i) \ g &\in H^{s}_{(\alpha , \beta)}(\mrI) \mbox{ then  } \ f g \in H^{s}_{(\alpha , \beta)}(\mrI), \ \mbox{ and  for } 
   \label{prodr1}  \\
(ii)  \ g &\in H^{-s}_{(\alpha , \beta)}(\mrI) \mbox{ then  } \ f g \in H^{-s}_{(\alpha , \beta)}(\mrI) . 
   \label{prodr2} 
\end{align}
\end{lemma}

\begin{theorem} \label{thmreg11}  \cite[Theorem 7.1]{erv191}
Let $s \ge -\alpha$, $\beta$ be determined by \textbf{Condition A},  \linebreak[4]
$c \in W^{\ceil{ \min \{s \, , \, \alpha \, + \, (\alpha - \beta) \, + \, 1 \, , \,  \alpha \, + \, \beta \, + \, 1\}}  , \infty}_{w}(\mrI)$ 
satisfying $c(x) \ge 0$
and
\begin{equation}
f \in H^{-\alpha/2}(\mrI) \cap H^{s}_{(\beta \, , \, \alpha - \beta)}(\mrI).
\label{def4f1}
\end{equation}  
Then there exists a unique solution $u(x) \ = \ (1 - x)^{\alpha - \beta} \, x^{\beta} \, \phi(x)$, with  \linebreak[4]
$\phi(x) \in 
H^{\alpha \, + \, \min \{s \, , \, \alpha \, + \, (\alpha - \beta) \, + \, 1 \, , \,  \alpha \, + \, \beta \, + \, 1\}}_{(\alpha - \beta \, , \, \beta)}(\mrI)$,
 to
\begin{equation}
\mcL_{r}^{\alpha} u(x) \ + \ c(x) \, u(x) \ = \ f(x) \, , \ x \in \mrI , \ 
\mbox{ subject to } u(0) = u(1) = 0 \, .
\label{theq1}
\end{equation}
\end{theorem}


The inclusion of an advection term can significantly reduced the regularity of the solution.

\begin{theorem} \label{thmreg13}  \cite[Theorem 7.2]{erv191}
Let $s \ge -\alpha$, $\beta$ be determined by \textbf{Condition A},  \linebreak[4]
$b , \, c \in W^{\ceil{\min \{s \, , \, \alpha \, + \, (\alpha - \beta) \, - \, 1 \, , \,  \alpha \, + \, \beta \, - \, 1\}} , \infty}_{w}(\mrI)$ satisfying 
$c(x)  \, - \, 1/2 D b(x) \ \ge 0$,
and
\begin{equation}
f \in H^{-\alpha/2}(\mrI) \cap H^{s}_{(\beta \, , \, \alpha - \beta)}(\mrI).
\label{def4f2}
\end{equation}  
Then there exists a unique solution $u(x) \ = \ (1 - x)^{\alpha - \beta} \, x^{\beta} \, \phi(x)$, with  \linebreak
$\phi(x) \in 
H^{\alpha \, + \, \min \{s \, , \, \alpha \, + \, (\alpha - \beta) \, - \, 1 \, , \,  \alpha \, + \, \beta \, - \, 1\}}_{(\alpha - \beta \, , \, \beta)}(\mrI)$,
to
\begin{equation}
\mcL_{r}^{\alpha} u(x) \ + \ b(x) \, D u(x) \ + \ c(x) \, u(x) \ = \ f(x) \, , \ x \in \mrI , \ 
\mbox{ subject to } u(0) = u(1) = 0 \, .
\label{theq2}
\end{equation}
\end{theorem}

Introduce $\wtilde{s}$ defined by
\be
\wtilde{s} \, := \, \left\{ \begin{array}{rl}
\min\{s , \, \alpha + (\alpha - \beta) +1 , \, \alpha +  \beta + 1 \} , \ & \mbox{ if } b = 0 \ \mbox{ (see Theorem \ref{thmreg11})} \\
\min\{s , \, \alpha + (\alpha - \beta) - 1 , \, \alpha +  \beta - 1 \} , \ & \mbox{ if } b \ne 0 \ \mbox{ (see Theorem \ref{thmreg13})} \, .
\end{array} \right.
\label{defsar}
\ee

 \setcounter{equation}{0}
\setcounter{figure}{0}
\setcounter{table}{0}
\setcounter{theorem}{0}
\setcounter{lemma}{0}
\setcounter{corollary}{0}
\setcounter{definition}{0}
%
\section{Weak Formulation}
\label{sec_form}
In place of \eqref{DefProb2}, \eqref{DefBC2}, we consider the following problem. 

\textit{Given $f \in H^{-\alpha/2}(\mrI) \cap H^{-\alpha/2}_{\omega^{*}}(\mrI)$,  and $b$ and $c$ satisfying the 
hypothesis of Theorem \ref{thmreg13},
determine $\phi \in H^{\alpha/2}_{\omega}(\mrI)$
such that $u(x) \ = \ \omega(x) \, \phi(x)$ satisfies}
\be
\langle \mcL_{r}^{\alpha} u \ + \ b \, D u \ +  \ c \, u \, , \, \psi \rangle_{\omega^{*}} 
\ = \ \langle f \, , \, \psi \rangle_{\omega^{*}} \, , \ \ \forall \, \psi \in H^{\alpha/2}_{\omega^{*}}(\mrI) \, .
\label{wform1}
\ee

Note that the formulation \eqref{wform1} has different test and trial spaces. With this in mind we recall the
Banach-Ne\v{c}as-Babu\v{s}ka theorem.
\begin{theorem} \cite[Pg. 85, Theorem 2.6]{ern041} \label{BNBthm}
Let $H_{1}$ and $H_{2}$ denote two real Hilbert spaces, $B(\cdot , \cdot) \, : \, H_{1} \times H_{2} \rightarrow \real$
a bilinear form, and $F \, : \, H_{2} \rightarrow \real$ a bounded linear functional on $H_{2}$. Suppose there are constants
$C_{1} < \infty$ and $C_{2} > 0$ such that
\begin{align}
&(i) \ | B(w , v) | \, \le \, C_{1} \, \| w \|_{H_{1}} \, \| v \|_{H_{2}} \, ,  \ 
\mbox{ for all } w \in H_{1} \, , \ v \in H_{2} \, , \label{BB1}  \\
&(ii) \ \sup_{0 \ne v \in H_{2}} \frac{ | B(w , v) |}{ \| v \|_{H_{2}} } \ \ge \ C_{2} \,  \| w \|_{H_{1}} \, ,  
 \mbox{ for all } w \in H_{1} \, , \label{BB2}  \\
&(iii) \ \sup_{w \in H_{1}}  | B(w , v) | \ > \ 0 \, ,  \mbox{ for all } v \in H_{2} \, , \ v \neq 0 \, .  \label{BB3}  
\end{align}
Then there exists a unique solution $w_{0} \in H_{1}$ satisfying $B(w_{0} \, , \, v) \, = \, F(v)$ for all 
$ v \in H_{2}$. Further, $\| w_{0} \|_{H_{1}} \, \le \, C_{2} \| F \|_{H_{2}}$.
\end{theorem}

For $f \in H^{-\alpha/2}(\mrI) \cap H^{-\alpha/2}_{\omega^{*}}(\mrI)$,  and $b$ and $c$ satisfying the 
hypothesis of Theorem \ref{thmreg13}, let $B \, : \, H^{\alpha/2}_{\omega} \times H^{\alpha/2}_{\omega^{*}} \rightarrow \real$, and
$F \, :  \, H^{\alpha/2}_{\omega^{*}} \rightarrow \real$ be defined by
\begin{align}
 B(\phi , \psi) &:= \ \langle \mcL_{r}^{\alpha} u \ + \ b \, D u \ +  \ c \, u \, , \, \psi \rangle_{\omega^{*}}  \, ,
 \label{defB}  \\
F(\psi) &:= \ \langle f \, , \, \psi \rangle_{\omega^{*}} \, .   \label{defF}
\end{align}

\subsection{Continuity of $B(\cdot , \cdot)$}
\label{ssec_CntB}
In order to establish that $B(\cdot , \cdot)$ is well defined and continuous we need to 
determine which $H^{t}_{(a , b)}(\mrI)$ space $ \omega \, \phi$ lies in.

\begin{figure}[!ht]
\begin{center}
 \includegraphics[height=2.5in]{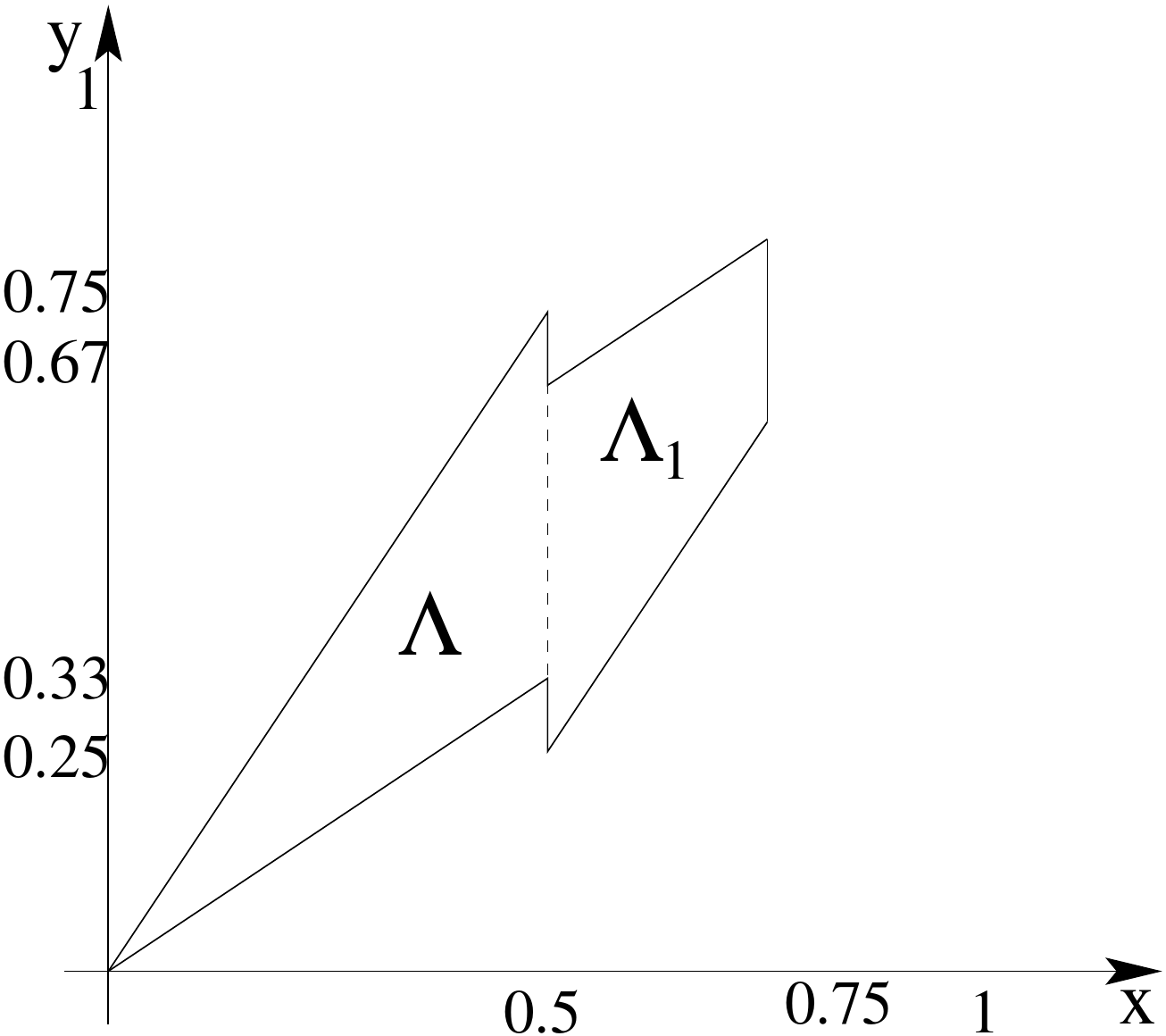}
   \caption{Domain $\Lambda^{*} \, = \, \Lambda \cup \Lambda_{1}$.}
   \label{figdomlams}
\end{center}
\end{figure}

The $H^{s}_{(a , b)}(\mrI)$ space a function $f$ lies in is determined by its behavior at: (i) the left endpoint ($x = 0$), 
(ii) the right endpoint ($x = 1$), and (iii) away from the endpoints. In order to separate the consideration of the endpoint 
behaviors, following \cite{ber921}, we introduce the following function space $H^{s}_{(\gamma)}(\mrJ)$. 
Let $\mrJ \, := \, (0 , \, 3/4)$, and
\begin{align*}
\Lambda^{*} &:= \ 
\left\{ (x , y) \, : \, \frac{2}{3} x < y < \frac{3}{2} x , \, 0 < x < \frac{1}{2} \right\} \cup
\left\{ (x , y) \, : \, \frac{3}{2} x - \frac{1}{2} \,  < y <  \, \frac{2}{3} x + \frac{1}{3} , \, 1/2 \le x < 3/4 \right\} \,  \\
&:= \ \Lambda \cup \Lambda_{1}  \ \ \mbox{(see Figure \ref{figdomlams})} \, .
\end{align*}

Introduce the semi-norm and norm
\begin{align}
 | f |_{H^{s}_{(\gamma)}(\mrJ)}^{2} &:= \ 
  \iint_{\Lambda}  x^{\gamma + s} \, 
  \frac{ | D^{\floor{s}} f(x) \, - \, D^{\floor{s}} f(y) |^{2}}{ | x \, - \, y |^{1 \, + \, 2 (s - \floor{s})}} dy \,  dx  \ + \ 
  \iint_{\Lambda_{1}}  x^{\gamma + s} \, 
  \frac{ | D^{\floor{s}} f(x) \, - \, D^{\floor{s}} f(y) |^{2}}{ | x \, - \, y |^{1 \, + \, 2 (s - \floor{s})}} dy \,  dx \, \nonumber \\
 &:=   | f |_{H^{s}_{(\gamma)}(\Lambda)}^{2}  \ + \  | f |_{H^{s}_{(\gamma)}(\Lambda_{1})}^{2}  \, , \nonumber \\
 &   \nonumber \\
\mbox{and } \ 
  \| f \|_{H^{s}_{(\gamma)}(\mrJ)}^{2} &:= \  \left\{ \begin{array}{rl}  
  \sum_{j = 0}^{s} \| D^{j} f \|_{L^{2}_{(\gamma + j)}(\mrJ)}^{2} \, , & \mbox{ for } s \in \mathbb{N}_{0}  \nonumber \\
 \sum_{j = 0}^{\floor{s}} \| D^{j} f \|_{L^{2}_{(\gamma + j)}(\mrJ)}^{2} \ + \  | f |_{H^{s}_{(\gamma)}(\mrJ)}^{2} \, , & \mbox{ for } s \in \mathbb{R}^{+} \backslash \mathbb{N}_{0} 
 \end{array} \right. \, ,   \label{deffnJ} \\
  &  \nonumber \\
\mbox{where } \ 
\| g \|_{L^{2}_{(\gamma)}(\mrJ)}^{2} &:= \ 
  \int_{\mrJ}  x^{\gamma} \, g^{2}(x) \, dx \, .  \nonumber
\end{align}

Then, $H^{s}_{(\gamma)}(\mrJ) \, := \, \{ f \, : \, f \mbox{ is measurable and }  \| f \|_{H^{s}_{(\gamma)}(\mrJ)} \, < \, \infty \}$. 

\textbf{Note}: A function $f(x)$ is in $H^{s}_{(a , b)}(\mrI)$ if and only if $f(\frac{3}{4} x) \in H^{s}_{(b)}(\mrJ)$ and
$f(\frac{3}{4}(1 - x)) \in H^{s}_{(a)}(\mrJ)$.

From \cite{erv191} we have the following theorem.
\begin{theorem} \cite[Theorem 6.4]{erv191} \label{thrmall1}
Let $n \, \le \,  s < n + 1$, $n \in \mathbb{N}_{0}$, $p \ge n$, $\mu > -1$, and $\psi \in H^{s}_{(\mu)}(\mrJ)$. 
Then $x^{p} \, \psi \in H^{t}_{(\sigma)}(\mrJ)$
provided
\be
 0 \le  t \le s \, , \ \ \sigma \, + \, 2 p \, \ge  \, \mu  \, , \  \ \sigma \, + \, 2 p \, - t \, >  \,  - 1 \, , \ 
   \ \mbox{ and } \ \ \sigma \, + \, 2 p \, + t \, \ge \, \mu \, + \, s   \, .
 \label{piu101}
\ee
Additionally, when \eqref{piu101} is satisfied, there exists $C > 0$ (independent of $\psi$) such that
$\| x^{p} \, \psi \|_{H^{t}_{(\sigma)}(\mrJ)} \, \le \, C \, \| \psi \|_{H^{s}_{(\mu)}(\mrJ)}$.
\end{theorem}

\begin{lemma} \label{lmct1}
The terms $\langle \mcL_{r}^{\alpha} \omega \phi  , \, \psi \rangle_{\omega^{*}}$,
$\langle  b \, D \omega \phi  \, , \, \psi \rangle_{\omega^{*}}$ and 
$\langle  c \, \omega \phi \, , \, \psi \rangle_{\omega^{*}}$ are well defined. Additionally, there exists $C > 0$ such that
for $\phi(x) \in H^{\alpha/2}_{\omega}(\mrI)$
and 
$\psi(x)  \in H^{\alpha/2}_{\omega^{*}}(\mrI)$
\be
|B( \phi , \psi )| \ = \ \left| \langle \mcL_{r}^{\alpha} \omega \phi \, + \, 
  b \, D \omega \phi   \, + \, 
   c \, \omega \phi \, , \, \psi \rangle_{\omega^{*}} \right| \ \le \ 
 C \, \| \phi \|_{H^{\alpha/2}_{\omega}(\mrI)}  \, \| \psi \|_{H^{\alpha/2}_{\omega^{*}}(\mrI)} \, . \label{cnty1}
 \ee
\end{lemma}
\textbf{Proof}: We begin by considering the $\langle  b \, D \omega \phi  \, , \, \psi \rangle_{\omega^{*}}$ term.

From Theorem \ref{thrmall1}, 
with $s \, = \, \alpha/2$, $\mu = \beta$, $p = \beta$, and choosing $\sigma \, = \, \alpha - \beta - 1$ we have that
$t \, \le \, \alpha/2$. Hence for $\phi_{0} \in H^{\alpha/2}_{\beta}(\mrJ)$, 
$x^{\beta} \phi_{0}(x) \in H^{\alpha/2}_{(\alpha - \beta - 1)}(\mrI)$, with 
$\| x^{\beta} \phi_{0}(x) \|_{H^{\alpha/2}_{(\alpha - \beta - 1)} }
\lesssim \| \phi_{0}(x) \|_{H^{\alpha/2}_{(\beta)}}$.

Again, using Theorem \ref{thrmall1}, 
with $s \, = \, \alpha/2$, $\mu \, = \, \alpha - \beta$, $p \, = \, \alpha - \beta$, 
and choosing $\sigma  =  \beta$ we have that
$t \, \le \, \alpha/2$. Hence for $\phi_{1} \in H^{\alpha/2}_{\alpha - \beta}(\mrJ)$, 
$x^{\alpha - \beta} \phi_{1}(x) \in H^{\alpha/2}_{(\beta - 1)}(\mrI)$ with 
$\| x^{\alpha - \beta} \phi_{1}(x) \|_{H^{\alpha/2}_{(\beta - 1)}} 
\lesssim \| \phi_{1}(x) \|_{H^{\alpha/2}_{(\alpha - \beta)}}$.

Combining the above two applications of Theorem \ref{thrmall1} we have that for $\phi \in H^{\alpha/2}_{\omega}(\mrI)$,
$\omega \phi \in H^{\alpha/2}_{(\beta - 1 \, , \, \alpha - \beta - 1)}(\mrI)$ with 
\be
  \| \omega \phi \|_{H^{\alpha/2}_{(\beta - 1 \, , \, \alpha - \beta - 1)}(\mrI)} 
  \ \lesssim \  \| \phi \|_{H^{\alpha/2}_{\omega}(\mrI)} \, .  
  \label{plok1}
\ee

A similar application of  Theorem \ref{thrmall1}  establishes that for $\phi \in H^{\alpha/2}_{\omega}(\mrI)$,
$\omega \phi \in H^{\alpha/2}_{\omega^{*}}(\mrI)$ with 
\be
  \| \omega \phi \|_{H^{\alpha/2}_{\omega^{*}}(\mrI)} 
  \ \lesssim \  \| \phi \|_{H^{\alpha/2}_{\omega}(\mrI)} \, .  
  \label{plok2}
\ee

From \eqref{plok1} and Lemma \ref{lmamapD}
we have that $D \omega \phi \in H^{\alpha/2 - 1}_{(\beta \, , \, \alpha - \beta)}(\mrI)$ with
$\|  D \omega \phi \|_{H^{\alpha/2 - 1}_{(\beta \, , \, \alpha - \beta)}} \lesssim \| \phi \|_{H^{\alpha/2}_{\omega}}$.
Thus, with the assumption on $b$ and using Lemma \ref{lmaprodsp},
\begin{align}
\langle  b \, D \omega \phi  \, , \, \psi \rangle_{\omega^{*}}
&\le \ \| D \omega \phi \|_{H^{\alpha/2 - 1}_{(\beta \, , \, \alpha - \beta)}} \, 
    \| b \, \psi \|_{H^{1 - \alpha/2}_{(\beta \, , \, \alpha - \beta)}}   \nonumber \\
&\lesssim \  \| \phi \|_{H^{\alpha/2}_{\omega}} \, \| \psi \|_{H^{1 - \alpha/2}_{\omega^{*}}} \, , \label{cnty3m1} \\
&\lesssim \  \| \phi \|_{H^{\alpha/2}_{\omega}} \, \| \psi \|_{H^{\alpha/2}_{\omega^{*}}} \, , \label{cnty3}
\end{align}
where in the last step we have used  $1 - \alpha/2 \, \le \, \alpha/2$.

For $\| \phi \| \in H^{\alpha/2}_{\omega}(\mrI)$ and $\| \psi \| \in H^{\alpha/2}_{\omega^{*}}(\mrI)$, using
\eqref{plok2} and the assumption on $c$,
\begin{align}
\langle  c \, \omega \phi \, , \, \psi \rangle_{\omega^{*}} 
&= \ \int_{I} \omega^{*}(x) \, c(x) \, \omega(x) \phi(x) \,  \psi(x) \, dx  \nonumber \\
&\le \ \|  \omega^{1/2} \,  \omega^{* 1/2} \, \|_{L^{\infty}}
 \int_{I}  \omega^{1/2}(x) \phi(x) \,  \omega^{* 1/2}(x) \, c(x) \, \psi(x) \, dx \nonumber  \\
&\le \  \| \phi \|_{L^{2}_{\omega}} \, \| c \, \psi \|_{L^{2}_{\omega^{*}}} \nonumber  \\
&\lesssim \| \phi \|_{H^{\alpha/2}_{\omega}} \, \| \psi \|_{H^{1 - \alpha/2}_{\omega^{*}}}   \label{cnty4m1}. \\
&\lesssim \| \phi \|_{H^{\alpha/2}_{\omega}} \, \| \psi \|_{H^{\alpha/2}_{\omega^{*}}} \, . \label{cnty4}
\end{align}

For $\phi(x) \, = \, \sum_{i = 0}^{\infty} \phi_{i} \, \what{G}_{i}^{(\alpha - \beta \, , \, \beta)}(x) \in H^{\alpha/2}_{\omega}(\mrI)$
and 
$\psi(x) \, = \, \sum_{j = 0}^{\infty} \psi_{j} \, \what{G}_{j}^{(\beta \, , \, \alpha - \beta)}(x) \in H^{\alpha/2}_{\omega^{*}}(\mrI)$,
using \eqref{propmcL}
\begin{align}
\langle \mcL_{r}^{\alpha} \omega \phi  , \, \psi \rangle_{\omega^{*}} 
&= \ \left( \sum_{i = 0}^{\infty} -c_{*}^{*} \, \lambda_{i} \, \phi_{i} \, \what{G}_{i}^{(\beta \, , \, \alpha - \beta)}(x) \, , \, 
\sum_{j = 0}^{\infty} \psi_{j} \, \what{G}_{j}^{(\beta \, , \, \alpha - \beta)}(x) \right)_{\omega^{*}}  \label{cnty45} \\
&= \ -c_{*}^{*} \sum_{k = 0}^{\infty} \lambda_{k} \, \phi_{k} \, \psi_{k}  \ 
\sim \ \sum_{k = 0}^{\infty} k^{\alpha} \, \phi_{k} \, \psi_{k}  \nonumber \\
&\lesssim \ \left( \sum_{k = 0}^{\infty} k^{\alpha} \, \phi_{k} \right)^{1/2} \, 
  \left( \sum_{k = 0}^{\infty} k^{\alpha} \, \psi_{k} \right)^{1/2}  
 \ \lesssim \ \left( \sum_{k = 0}^{\infty} (1 + k^{2})^{\alpha/2} \, \phi_{k} \right)^{1/2} \, 
  \left( \sum_{k = 0}^{\infty} (1 + k^{2})^{\alpha/2} \, \psi_{k} \right)^{1/2}  \nonumber  \\
&\lesssim \ \| \phi \|_{H^{\alpha/2}_{\omega}} \, \| \psi \|_{H^{\alpha/2}_{\omega^{*}}} \, , 
\ \mbox{ using } \eqref{defHr} \, .  \label{cnty5}
\end{align}

Combining \eqref{cnty3}, \eqref{cnty4} and \eqref{cnty5} we obtain \eqref{cnty1}. \\
\mbox{ } \hfill \qed

\subsection{Conditions \eqref{BB2} and \eqref{BB3}}
\label{ssec_BCoe}
For the case $r = 1/2$ we have $\alpha - \beta \, = \, \beta \, = \, \alpha/2$ and, consequently,
$\omega = \omega^{*}$. In this case for $\psi = \phi$
\begin{align}
\langle b \, D (\omega \, \phi) \ + \ c \, \omega \, \phi \ , \ \psi \rangle_{\omega}
&= \ \int_{0}^{1} \omega \left( b \, D (\omega \, \phi) \ + \ c \, \omega \, \phi \right) \phi \ dx  \nonumber  \\
&= \ \int_{0}^{1} b \, \frac{1}{2} D (\omega \, \phi)^{2} \ + \ c \, ( \omega \, \phi )^{2} \ dx   \nonumber  \\
&= \ \int_{0}^{1} \left( c \, - \, \frac{1}{2} D b \right) \ ( \omega \, \phi )^{2} \ dx  \, . \label{coe1}
\end{align}
Proceeding as in \eqref{cnty45}, for $\psi = \phi$ and $\omega^{*} = \omega$,
\begin{align}
\langle \mcL_{1/2}^{\alpha} (\omega \, \phi)  , \, \phi \rangle_{\omega} 
&\sim \ \sum_{k = 0}^{\infty} k^{\alpha} \, \phi_{k}^{2} 
\ \sim \ \sum_{k = 0}^{\infty} (1 + k^{2})^{\alpha/2} \, \phi_{k}^{2} \nonumber \\
&\sim \ \| \phi \|^{2}_{H^{\alpha/2}_{(\alpha/2 , \alpha/2)}} \, .   \label{coe2}
\end{align}
Hence for $( c \, - \, \frac{1}{2} D b ) \ge 0$, combining \eqref{coe1} and \eqref{coe2} we have that
$B(\cdot , \cdot)$ is coercive on $H^{\alpha/2}_{(\alpha/2 , \alpha/2)} \times H^{\alpha/2}_{(\alpha/2 , \alpha/2)}$.
Then, from the Lax-Milgram, we have the following lemma.
\begin{lemma} \label{exreq12}
For $1 < \alpha < 2$ and $r = 1/2$, given $f \in H^{-\alpha/2}_{(\alpha/2 , \alpha/2)}(\mrI)$ and
$b(x)$ and $c(x)$ satisfying $c(x) \ - \ 1/2 \, D b(x) \, \ge \, 0 , \ x \in \mrI$, there exists a unique solution
$u(x) \ = \ (1 - x)^{\alpha/2} x^{\alpha/2} \, \phi(x)$ to \eqref{DefProb2}, \eqref{DefBC2}, with
$\phi \in H^{\alpha/2}_{(\alpha/2 , \alpha/2)} (\mrI)$ satisfying 
$\| \phi \|_{H^{\alpha/2}_{(\alpha/2 , \alpha/2)}(\mrI)} \lesssim \| f \|_{H^{-\alpha/2}_{(\alpha/2 , \alpha/2)}(\mrI)}$. 
\end{lemma}

This special case of \eqref{wform1} corresponding to $r = 1/2$ has been thoroughly investigated by Hao and Zhang in \cite{hao201}.

For the general case, $(r \neq \frac{1}{2})$, to show \eqref{BB2} and \eqref{BB3}, and hence establish the
well posedness of the formulation, following an approach by Jin, Lazarov and Zhou in \cite{jin162}, we use
the Petree-Tartar Lemma.

\begin{lemma} \label{lemPT} \cite[Pg. 469]{ern041} (Petree-Tartar). Let $X$, $Y$, $Z$ be three Banach spaces. Let
$A \in \mcL( X ; Y)$ be an injective operator and let $T \in \mcL(X ; Z)$ be a compact operator. If there
exists $c_{1} > 0$ such that $c_{1} \, \| x \|_{X} \, \le \, \| A x \|_{Y} \, + \, \| T x \|_{Z}$, then $\mathrm{Im}(A)$ is
closed; equivalently, there is $c_{2} > 0$ such that
\be
  \forall x \in X , \ \ \ c_{2} \| x \|_{X} \, \le \, \| A x \|_{Y} \, .
  \label{PTlem}
\ee
\end{lemma}  

To relate the Petree-Tartar Lemma to the formulation \eqref{wform1}, let $X = H^{\alpha/2}_{\omega}(\mrI)$,
$Y = Z = H^{-\alpha/2}_{\omega^{*}}(\mrI)$,
\begin{align*}
&A  \, : \, X \rightarrow Y \ \ \mbox{ be defined by } \ \ 
A \phi \ := \ \mcL_{r}^{\alpha} \omega \phi \ + \ b \, D  \omega \phi \ + \ c \,  \omega \phi \, ,
\ \ \mbox{ and }   \\
&T  \, : \, X \rightarrow Z \ \ \mbox{ be defined by } \ \ T \phi \ := \ - \left( \ b \, D  \omega \phi \ + \ c \,  \omega \phi  \right) \, .
\end{align*}
That $A \in \mcL( X ; Y)$ follows from its definition and the continuity of $B(\cdot , \cdot)$. Its injectivity follows from the
uniqueness of solution to \eqref{theq2}. 
The fact that $T \in \mcL( X ; Z)$ follows from its definition and \eqref{cnty3} and \eqref{cnty4}. Also, from 
\eqref{cnty3m1} and \eqref{cnty4m1} we have that 
$T \, : \, H^{\alpha/2}_{\omega}(\mrI) \rightarrow H^{1 - \alpha/2}_{\omega^{*}}(\mrI)$ is bounded. As 
$H^{s}_{\omega^{*}}(\mrI)$ is compactly embedded in $H^{t}_{\omega^{*}}(\mrI)$ for $s > t$, 
\cite[pg. 10, Remark 2]{erv191}, since
$1 \, - \, \alpha/2 \, > \, - \alpha/2$, it follows that $T \in \mcL(X ; Z)$ is a compact operator.

Let $\phi(x) \ = \ \sum_{i = 1}^{\infty} \phi_{i} \what{G}_{i}^{(\alpha - \beta , \beta)}(x) \in H^{\alpha/2}_{\omega}(\mrI)$ and
$\psi(x) \ = \ \sum_{i = 1}^{\infty} \phi_{i} \what{G}_{i}^{(\beta , \alpha - \beta)}(x) \in H^{\alpha/2}_{\omega^{*}}(\mrI)$.
Note that $\| \phi \|_{H^{\alpha/2}_{\omega}(\mrI)} \, = \, \| \psi \|_{H^{\alpha/2}_{\omega^{*}}(\mrI)}$. Then,
\begin{align}
\| \phi \|_{H^{\alpha/2}_{\omega}(\mrI)}^{2} &= \ \sum_{i = 0}^{\infty} \big( 1 + i^{2} \big)^{\alpha/2} \phi_{i}^{2}   \nonumber \\
&\lesssim \ \sum_{i = 0}^{\infty} \lambda_{i} \,  \phi_{i}^{2} \ = \ \langle \mcL_{r}^{\alpha} \omega \phi \, , \, \psi \rangle_{\omega^{*}}
 \nonumber \\
&= \  \langle \mcL_{r}^{\alpha} \omega \phi \ + \ b \, D  \omega \phi \ + \ c \,  \omega \phi \, , \, \psi \rangle_{\omega^{*}}
 \ + \  \langle  - \left( \ b \, D  \omega \phi \ + \ c \,  \omega \phi \right) \, , \, \psi \rangle_{\omega^{*}}. \nonumber \\
&= \   \langle A \phi \, , \, \psi \rangle_{\omega^{*}} \ + \ \langle T \phi \, , \, \psi \rangle_{\omega^{*}}  \nonumber \\      
&\le  \| A \phi \|_{H^{-\alpha/2}_{\omega^{*}}(\mrI)} \, \| \psi \|_{H^{\alpha/2}_{\omega^{*}}(\mrI)} \ + \ 
   \| T \phi \|_{H^{-\alpha/2}_{\omega^{*}}(\mrI)} \, \| \psi \|_{H^{\alpha/2}_{\omega^{*}}(\mrI)} \, .  \nonumber 
\end{align}
Using $\| \phi \|_{H^{\alpha/2}_{\omega}(\mrI)} \, = \, \| \psi \|_{H^{\alpha/2}_{\omega^{*}}(\mrI)}$, we obtain that there exists $c_{1} > 0$
such that
\[
 c_{1} \,  \| \phi \|_{X} \, \le \, \, \| A \phi \|_{Y} \, + \, \| T \phi \|_{Z} \, .
\] 
Then, applying the Petree-Tartar Lemma, it follows that there exists $C_{2} > 0$ such that
\be
C_{2} \| \phi \|_{X} \, \le \, \| A \phi \|_{Y} \, , \ \  \mbox{i.e.}, \ \ 
C_{2} \| \phi \|_{H^{\alpha/2}_{\omega}(\mrI)} \, \le \, \| A \phi \|_{H^{- \alpha/2}_{\omega^{*}}(\mrI)} \, .
\label{PTn1}
\ee

\begin{lemma} \label{Bcondii}
For $B(\cdot , \cdot)$ defined by \eqref{defB}, the
condition (ii) given by \eqref{BB2} is satisfied.
\end{lemma}
\textbf{Proof}: Noting that
\[
\sup_{0 \ne v \in H^{\alpha/2}_{\omega^{*}}(\mrI)}
\frac{ | B(w , v) |}{ \| v \|_{H^{\alpha/2}_{\omega^{*}}(\mrI)} } \ \ge \ C_{2} \,  \| w \|_{H^{\alpha/2}_{\omega}(\mrI)}  \ \ 
\mbox{ is equivalent to } \ \
\| A w \|_{H^{- \alpha/2}_{\omega^{*}}(\mrI)}  \ \ge \ C_{2} \,  \| w \|_{H^{\alpha/2}_{\omega}(\mrI)} \, ,
\]
the condition (ii) follows from \eqref{PTn1}.  \\
\mbox{ } \hfill \qed

\begin{lemma} \label{Bcondiii}
For $B(\cdot , \cdot)$ defined by \eqref{defB}, the
condition (iii) given by \eqref{BB3} is satisfied.
\end{lemma}
\textbf{Proof}: 
The adjoint problem to \eqref{wform1} is: 
\textit{Given $g \in H^{-\alpha/2}_{\omega}(\mrI)$, determine $\psi \in H^{\alpha/2}_{\omega^{*}}(\mrI)$
such that $v(x) \ = \ \omega^{*}(x) \, \psi(x)$ satisfies}
\be
\langle \mcL_{(1 - r)}^{\alpha} v \ - \ b \, D v \ +  \ (c \, - \, D b) v \, , \, \phi \rangle_{\omega} 
\ = \ \langle g \, , \, \phi \rangle_{\omega} \, , \ \ \forall \, \phi \in H^{\alpha/2}_{\omega}(\mrI) \, .
\label{wformc}
\ee
This weak form corresponds to the fractional diffusion, advection, reaction equation (see Theorem \ref{thmreg13}):
\textit{Given $g \in H^{-\alpha/2}(\mrI) \cap H^{-\alpha/2}_{\omega}(\mrI)$ determine $v(x)$ satisfying}
\begin{equation}
\mcL_{(1 - r)}^{\alpha} v(x) \ - \ b(x) \, D v(x) \ + \ \big( c(x) \, - \, D b(x) \big) v(x) \ = \ g(x) \, , \ x \in \mrI , \ 
\mbox{ subject to } v(0) = v(1) = 0 \, .
\label{theq2adj}
\end{equation}
An analogous argument as used to establish condition (ii) given by \eqref{BB2} can be applied to 
\eqref{wformc} and \eqref{theq2adj} to establish condition (iii) given by \eqref{BB3}. \\
\mbox{ } \hfill \qed

Combining Lemmas \ref{lmct1}, \ref{Bcondii} and \ref{Bcondiii} with Theorem \ref{BNBthm} we obtain the following.
\begin{theorem} \label{thrmexw1}
There exists a unique solution $\phi$ to \eqref{wform1}, satisfying 
$\|\phi \|_{H^{\alpha/2}_{\omega}(\mrI)} \, \le \, \frac{1}{C_{2}} \| f \|_{H^{- \alpha/2}_{\omega^{*}}(\mrI)}$.
\end{theorem}
\textbf{Proof}: First, note that $F$ defined by \eqref{defF} satisfies
\begin{align*}
\| F \| &= \ \sup_{0 \ne \psi \in H^{\alpha/2}_{\omega^{*}}(\mrI)} \frac{ | F(\psi) |}{ \| \psi \|_{H^{\alpha/2}_{\omega^{*}}(\mrI)}}
\ = \ \sup_{0 \ne \psi \in H^{\alpha/2}_{\omega^{*}}(\mrI)}
 \frac{ | \langle f , \psi \rangle_{\omega^{*}} |}{ \| \psi \|_{H^{\alpha/2}_{\omega^{*}}(\mrI)}}    \\
&\le \ \sup_{0 \ne \psi \in H^{\alpha/2}_{\omega^{*}}(\mrI)}
 \frac{ \| f \|_{H^{\alpha/2}_{\omega^{*}}(\mrI)} \| \psi \|_{H^{\alpha/2}_{\omega^{*}}(\mrI)} }%
 { \| \psi \|_{H^{\alpha/2}_{\omega^{*}}(\mrI)}}
 \ = \ \| f \|_{H^{\alpha/2}_{\omega^{*}}(\mrI)} \, .
 \end{align*}
Hence, $F$ defines a bounded linear functional.
The existence and uniqueness of $\phi$ then follows from 
combining Lemmas \ref{lmct1}, \ref{Bcondii} and \ref{Bcondiii} with Theorem \ref{BNBthm}. To obtain the bound
for $\|\phi \|_{H^{\alpha/2}_{\omega}(\mrI)}$, from Lemma \ref{Bcondii}
\begin{align}
\| \phi \|_{H^{\alpha/2}_{\omega}(\mrI)} 
&\le \ \frac{1}{C_{2}} \sup_{0 \ne \psi \in H^{\alpha/2}_{\omega^{*}}}
\frac{ | B(\phi , \psi_) |}{ \| \psi \|_{H^{\alpha/2}_{\omega^{*}}(\mrI)}} \ 
\ = \  \frac{1}{C_{2}}
\sup_{0 \ne \psi \in H^{\alpha/2}_{\omega^{*}}}
\frac{ | \langle f \, , \, \psi \rangle_{\omega^{*}} |}{ \| \psi \|_{H^{\alpha/2}_{\omega^{*}}(\mrI)}}  \nonumber \\
&\le \ \frac{1}{C_{2}} 
\sup_{0 \ne \psi \in H^{\alpha/2}_{\omega^{*}}}
\frac{ \| f \|_{H^{- \alpha/2}_{\omega^{*}}(\mrI)} \, \| \psi \|_{H^{\alpha/2}_{\omega^{*}}(\mrI)} }%
{ \| \psi \|_{H^{\alpha/2}_{\omega^{*}}(\mrI)}} 
\ = \ \frac{1}{C_{2}}  \| f \|_{H^{- \alpha/2}_{\omega^{*}}(\mrI)}  \, .  \label{esthg1}
\end{align}
\mbox{ } \hfill \qed

\begin{corollary}  \label{regH}
For $f \in H^{-\alpha/2}(\mrI) \cap H^{s}_{\omega^{*}}(\mrI)$,  $s \ge -\alpha/2$, and $b$ and $c$ satisfying the 
hypothesis of Theorem \ref{thmreg13},
there
exists $C > 0$ such that with $\phi$
 given by \eqref{wform1} satisfies
\be
 \| \phi  \|_{H^{\wtilde{s} +  \alpha}_{\omega}(\mrI)} \ \le \ C \, \| f \|_{H^{\wtilde{s}}_{\omega^{*}}(\mrI)} \, .
\label{hphibd}
\ee
\end{corollary}
\textbf{Proof}:
The proof follows as that for Theorems \ref{thmreg11} and  \ref{thmreg13}.
At each of the (finite number of) steps in the boot strapping argument
the terms on the right hand side are bounded by a constant times $ \| f \|_{H^{\wtilde{s}}_{\omega}(\mrI)}$. \\
\mbox{ } \hfill \qed \\

 \setcounter{equation}{0}
\setcounter{figure}{0}
\setcounter{table}{0}
\setcounter{theorem}{0}
\setcounter{lemma}{0}
\setcounter{corollary}{0}
\setcounter{definition}{0}
%
\section{Approximation Scheme}
\label{sec_appx}
As $ \{ \what{G}_{j}^{(a , b)} \}_{j = 0}^{\infty}$ is a basis for $H^{\alpha/2}_{(a , b)}(\mrI)$, let
$X_{N} \, := \, \mbox{span}  \{ \what{G}_{j}^{(\alpha - \beta \, , \, \beta)} \}_{j = 0}^{N} 
\subset H^{\alpha/2}_{(\alpha - \beta \, , \, \beta)}(\mrI)$, and 
$Y_{N} \, := \, \mbox{span}  \{ \what{G}_{j}^{(\beta \, , \, \alpha - \beta)} \}_{j = 0}^{N}
\subset H^{\alpha/2}_{(\beta \, , \, \alpha - \beta)}(\mrI)$. Corresponding to \eqref{wform1} we have 
the following approximation scheme.

\textit{Given $f \in H^{-\alpha/2}(\mrI) \cap H^{-\alpha/2}_{\omega^{*}}(\mrI)$,  and $b$ and $c$ satisfying the 
hypothesis of Theorem \ref{thmreg13}, determine $\phi_{N} \in X_{N}$
such that $u_{N}(x) \ = \ \omega(x) \, \phi_{N}(x)$ satisfies}
\be
\langle \mcL_{r}^{\alpha} \omega(x) \, \phi_{N}(x)  +  \ b(x) \, D \omega(x) \, \phi_{N}(x)  + 
 \ c(x) \, \omega(x) \, \phi_{N}(x) \, , \, \psi_{N} \rangle_{\omega^{*}} 
\ = \ \langle f \, , \, \psi_{N} \rangle_{\omega^{*}} \, , \ \ \forall \, \psi_{N} \in Y_{N} \, .
\label{appx1}
\ee

The following lemma is used to establish the well posedness of \eqref{appx1}.

\begin{lemma} \label{lmadisB}
There exists $C_{3} > 0$, such that for $N$ sufficiently large,
\begin{align}
\sup_{0 \ne \psi_{N} \in Y_{N}} 
\frac{ | B(\phi_{N} , \psi_{N}) |}{ \| \psi_{N} \|_{H^{\alpha/2}_{\omega^{*}}(\mrI)}} \ 
&\ge \ C_{3} \, \| \phi_{N} \|_{H^{\alpha/2}_{\omega}(\mrI)} \, , \ \ \forall \, \phi_{N} \in X_{N} \, . \label{ddrtgf1}  
\end{align}
\end{lemma}
\textbf{Proof}: Let $\phi_{N} \in X_{N}$. For $\psi \in H^{\alpha/2}_{\omega^{*}}(\mrI)$, let
$\psi_{N} \ = \ \sum_{i = 0}^{N} \psi_{i} G^{(\beta , \alpha - \beta)}_{i}(x)$. Using Lemma \ref{Bcondii},
\begin{align}
& C_{2} \,  \| \phi_{N} \|_{H^{\alpha/2}_{\omega}(\mrI)} 
\ \le  \ \sup_{0 \ne \psi \in H^{\alpha/2}_{\omega^{*}}(\mrI)}
\frac{ B(\phi_{N} , \psi) }{ \| \psi \|_{H^{\alpha/2}_{\omega^{*}}(\mrI)} }  
\ \le \ 
\sup_{0 \ne \psi \in H^{\alpha/2}_{\omega^{*}}(\mrI)}
\frac{  B(\phi_{N} , \psi_{N})  \ + \ B(\phi_{N} \, , \, \psi - \psi_{N}) }{ \| \psi \|_{H^{\alpha/2}_{\omega^{*}}(\mrI)} }    \nonumber \\
& \quad \le \ 
\sup_{0 \ne \psi \in H^{\alpha/2}_{\omega^{*}}(\mrI)}
\frac{  B(\phi_{N} , \psi_{N})  }{ \| \psi \|_{H^{\alpha/2}_{\omega^{*}}(\mrI)} }  
\ + \ 
\sup_{0 \ne \psi \in H^{\alpha/2}_{\omega^{*}}(\mrI)}
\frac{  \langle - \big(  b(x) \, D \omega(x) \, \phi_{N}(x)  + 
 \ c(x) \, \omega(x) \, \phi_{N}(x) \big) \, , \psi - \psi_{N} \rangle_{\omega^{*}} }{ \| \psi \|_{H^{\alpha/2}_{\omega^{*}}(\mrI)} }  \, ,  
 \label{thre3}
\end{align}
where in the last step we have used
$\langle \mcL_{r}^{\alpha} \phi_{N} \, , \, \psi - \psi_{N} \rangle_{\omega^{*}} \, = \, 0$. 

From \eqref{cnty3m1} and \eqref{cnty4m1}, and using \eqref{Approx},
\begin{align}
\big| \langle \big(  b(x) \, D \omega(x) \, \phi_{N}(x)  + 
 \ c(x) \, \omega(x) \, \phi_{N}(x) \big) \, , \psi - \psi_{N} \rangle_{\omega^{*}} \big|
 &\le \ C \,   \| \phi_{N} \|_{H^{\alpha/2}_{\omega}} \,  \| \psi - \psi_{N} \|_{H^{1 - \alpha/2}_{\omega^{*}}}  \nonumber \\
&\le \ C \,  \| \phi_{N} \|_{H^{\alpha/2}_{\omega}} \, N^{1 - \alpha} \, \| \psi \|_{H^{\alpha/2}_{\omega^{*}}} \, . 
  \label{thre4}
\end{align}

Combining \eqref{thre3} and \eqref{thre4}, for $N$ sufficiently large we obtain \eqref{ddrtgf1}. \\
\mbox{ } \hfill \qed

\begin{theorem} \label{exds1}
There exists a unique $\phi_{N} \in H^{\alpha/2}_{\omega}(\mrI)$ satisfying \eqref{appx1}. In addition, for $C_{3}$ given
in \eqref{ddrtgf1}, 
$\| \phi_{N} \|_{H^{\alpha/2}_{\omega}(\mrI)} \, \le \, \frac{1}{C_{3}} \| f \|_{H^{-\alpha/2}_{\omega^{*}}(\mrI)}$.
\end{theorem}
\textbf{Proof}: For $\phi_{N} \ = \ \sum{j = 0}^{N} c_{j} \, \what{G}^{(\alpha - \beta , \beta)}_{j}(x)$, from \eqref{appx1},
the constants $c_{j}$ are determined from
\[
    \mathbb{A} \bfc \, = \, \bfb \, , \ \mbox{ where } \ \mathbb{A}_{i + 1 \, j + 1} \, = \, 
    B( \what{G}^{(\alpha - \beta , \beta)}_{j} \, , \,  \what{G}^{(\beta , \alpha - \beta)}_{i} ) \, , \
    \mbox{ and } \bfb_{i} \, = \, \langle f(x) , \what{G}^{(\beta , \alpha - \beta)}_{i}(x) \rangle_{\omega^{*}} , \ 
\]    
for $0 \le i , j \le N$.
Condition \eqref{ddrtgf1} implies the invertible of the square matrix $ \mathbb{A}$, and hence the uniqueness of $\phi_{N}$
satisfying \eqref{appx1}.  The bound for $\phi_{N}$ is obtained in an analogous manner to the bound for $\phi$ in \eqref{esthg1}. \\
\mbox{ } \hfill \qed

For $\phi_{N}$ given by \eqref{appx1} we have the following error bound.
\begin{lemma} \label{errbd1}
There exists $C > 0$ such that for $\phi$ satisfying
\eqref{wform1} and $\phi_{N}$ satisfying \eqref{appx1}
\be
  \| \phi \, - \, \phi_{N} \|_{H^{\alpha/2}_{\omega}(\mrI)} \ \le \ C \ \inf_{\zeta_{N} \in X_{N}}
  \| \phi \, - \, \zeta_{N} \|_{H^{\alpha/2}_{\omega}(\mrI)} \, .
 \label{erest1}
\ee
\end{lemma}
\textbf{Proof}: Note that for $\zeta_{N} \in X_{N}$, using \eqref{ddrtgf1},
\begin{align}
C_{3} \,\| \phi_{N} \, - \, \zeta_{N} \|_{H^{\alpha/2}_{\omega}(\mrI)}  
&\le \  \sup_{\stackrel{\psi_{N} \in Y_{N}}{\psi_{N} \neq 0}} 
\frac{ | B(\phi_{N} \, - \, \zeta_{N} \,  , \, \psi_{N}) |}{ \| \psi_{N} \|_{H^{\alpha/2}_{\omega^{*}}(\mrI)}}
\ = \ \sup_{\stackrel{\psi_{N} \in Y_{N}}{\psi_{N} \neq 0}} 
\frac{ | \langle f \, , \, \psi_{N} \rangle_{\omega^{*}} \ - \ 
B(\zeta_{N} \,  , \, \psi_{N}) |}{ \| \psi_{N} \|_{H^{\alpha/2}_{\omega^{*}}(\mrI)}}  \nonumber \\
&= \  \sup_{\stackrel{\psi_{N} \in Y_{N}}{\psi_{N} \neq 0}} 
\frac{ | B(\phi \, - \, \zeta_{N} \,  , \, \psi_{N}) |}{ \| \psi_{N} \|_{H^{\alpha/2}_{\omega^{*}}(\mrI)}}
 \ \ \ \mbox{(using \eqref{wform1})}  \nonumber \\
&\le \sup_{\stackrel{\psi_{N} \in Y_{N}}{\psi_{N} \neq 0}} 
\frac{ C_{1} \, \| \phi \, - \, \zeta_{N}\|_{H^{\alpha/2}_{\omega}(\mrI)} \, 
 \| \psi_{N} \|_{H^{\alpha/2}_{\omega^{*}}(\mrI)}}{ \| \psi_{N} \|_{H^{\alpha/2}_{\omega^{*}}(\mrI)}}  
\ = \  C_{1} \, \| \phi \, - \, \zeta_{N}\|_{H^{\alpha/2}_{\omega}(\mrI)} \, .  \label{erest2}
\end{align}
With the triangle inequality and \eqref{erest2}, we obtain
\[
  \| \phi \, - \, \phi_{N} \|_{H^{\alpha/2}_{\omega}(\mrI)}
  \ \le \ \| \phi \, - \, \zeta_{N} \|_{H^{\alpha/2}_{\omega}(\mrI)} \ + \ 
    \| \zeta_{N} \, - \, \phi_{N} \|_{H^{\alpha/2}_{\omega}(\mrI)}
\ \le \ (1 \, + \, C_{1}) \,  \| \phi \, - \, \zeta_{N} \|_{H^{\alpha/2}_{\omega}(\mrI)} \, .    
\]
As $\zeta_{N} \in X_{N}$ is arbitrary, then \eqref{erest1} follows. \\
\mbox{ } \hfill \qed

Combining Lemma \ref{errbd1} with Lemma \ref{lem:Approx} and Theorems \ref{thmreg11} and \ref{thmreg13} we obtain the
following error estimate.
\begin{corollary}  \label{errH}
For $f \in H^{-\alpha/2}(\mrI) \cap H^{s}_{\omega^{*}}(\mrI)$,  $s \ge -\alpha/2$, and $b$ and $c$ satisfying the 
hypothesis of Theorem \ref{thmreg13},
there
exists $C > 0$ such that for $\phi$ satisfying
\eqref{wform1} and $\phi_{N}$ satisfying \eqref{appx1}
\be
 \| \phi \, - \, \phi_{N} \|_{H^{\alpha/2}_{\omega}(\mrI)} \ \le \ C \, N^{\wtilde{s} \, + \, \alpha/2} \, 
 \| \phi \|_{H^{\wtilde{s} \, + \, \alpha}_{\omega}(\mrI)} \,
  \le \, C \, N^{\wtilde{s} \, + \, \alpha/2} \, \| f \|_{H^{\wtilde{s}}_{\omega^{*}}(\mrI)}   \,  .
\label{herrest}
\ee
\end{corollary}
\textbf{Proof}: From Corollary \ref{regH} we have that $\phi$ satisfies
$\phi \in H^{\wtilde{s} \, + \, \alpha}_{\omega}(\mrI)$. Then,  
applying Lemma \ref{lem:Approx}, with $\mu \, = \, \alpha/2$ and $t \, = \, \wtilde{s} \, + \, \alpha$, and using  
Corollary \ref{regH}, we obtain
\eqref{herrest}.  \\
\mbox{ } \hfill \qed

An estimate for $\| \phi \, - \, \phi_{N} \|_{L^{2}_{\omega}(\mrI)}$ can be obtained using a Aubin-Nitsche type argument.
\begin{corollary}  \label{errL2}
For $f \in H^{-\alpha/2}(\mrI) \cap H^{s}_{\omega^{*}}(\mrI)$,  $s \ge -\alpha/2$, and $b$ and $c$ satisfying the 
hypothesis of Theorem \ref{thmreg13},
there
exists $C > 0$ such that for $\phi$ satisfying
\eqref{wform1} and $\phi_{N}$ satisfying \eqref{appx1}
\be
 \| \phi \, - \, \phi_{N} \|_{L^{2}_{\omega}(\mrI)} \ \le \ C \, N^{\wtilde{s} \, + \, \alpha} \, 
 \| \phi \|_{H^{\wtilde{s} \, + \, \alpha}_{\omega}(\mrI)} \,
  \le \, C \, N^{\wtilde{s} \, + \, \alpha} \, \| f \|_{H^{\wtilde{s}}_{\omega^{*}}(\mrI)} \,  .
\label{L2errest}
\ee
\end{corollary}
\textbf{Proof}:  Introduce $\psi \in H^{\alpha/2}_{\omega^{*}}(\mrI)$ satisfying
\[
\mcL^{\alpha}_{(1 - r)} \omega^{*} \, \psi \  -  \ b \, D \omega^{*} \, \psi \ + \ \big(c \, - \, D b) \, \omega^{*} \, \psi \ = \ \phi \, - \, \phi_{N} \, .
\]
As $(\phi \, - \, \phi_{N}) \in L^{2}_{\omega}(\mrI)$, analogous to \eqref{hphibd}, we have that
\be
  \| \psi \|_{H^{\alpha}_{\omega^{*}}(\mrI)} \ \le \ C \, \| \phi \, - \, \phi_{N} \|_{L^{2}_{\omega}(\mrI)} \, .
\label{ghre1}
\ee
Then,
\begin{align*}
 \| \phi \, - \, \phi_{N} \|_{L^{2}_{\omega}} 
 &= \ \big( (\phi \, - \, \phi_{N}) \, , \, (\phi \, - \, \phi_{N}) \big)_{\omega} 
 \ = \ \big( (\phi \, - \, \phi_{N}) \, , \, \mcL^{\alpha}_{(1 - r)} \omega^{*} \, \psi \  -  \ b \, D \omega^{*} \, \psi 
 \ + \ \big(c \, - \, D b) \, \omega^{*} \, \psi  \big)_{\omega}  \\
&= \  \big( \mcL^{\alpha}_{r} \omega \, (\phi \, - \, \phi_{N}) \ + \ b \, D \omega \, (\phi \, - \, \phi_{N}) \ 
+ \ c \, \omega \,  (\phi \, - \, \phi_{N})
\, , \,  \psi  \big)_{\omega^{*}}  \ = \ B((\phi \, - \, \phi_{N}) \, , \, \psi)  \\
&= \ B((\phi \, - \, \phi_{N}) \, , \, \psi \, - \, \eta_{N}) \, , \ \ \mbox{for } \eta_{N} \in Y_{N} \, , \ \mbox{(using Galerkin orthogonality)} \\
&\le \  C_{1} \, \| \phi \, - \, \phi_{N} \|_{H^{\alpha/2}_{\omega}} \, \| \psi \, - \, \eta_{N} \|_{H^{\alpha/2}_{\omega^{*}}} \, , 
\ \ \mbox{using } \eqref{cnty1} \, , \\
&\le \ C \, N^{\wtilde{s} \, + \, \alpha/2} \, 
 \| \phi \|_{H^{\wtilde{s} \, + \, \alpha}_{\omega}} \, N^{\alpha/2} \, \| \psi \|_{H^{\alpha}_{\omega^{*}}} \, ,
 \ \ \mbox{using } \eqref{Approx} \, ,  \\
&\le \ C \, N^{\wtilde{s} \, + \, \alpha} \, 
 \| \phi \|_{H^{\wtilde{s} \, + \, \alpha}_{\omega}} \,  \| \phi \, - \, \phi_{N} \|_{L^{2}_{\omega}}  \, ,
 \ \ \mbox{using } \eqref{ghre1} \, .
\end{align*}
Finally, dividing through by $\| \phi \, - \, \phi_{N} \|_{L^{2}_{\omega}}$ and using \eqref{hphibd} we obtain \eqref{L2errest}. \\
\mbox{ } \hfill \qed

\textbf{Error estimate for $u - u_{N}$}. \\
The weighted $L^{2}_{\omega^{-1}}$ error estimate for $u - u_{N}$, where $u_{N} \, := \, \omega \, \phi_{N}$, follows easily
from the definitions of $u_{N}$ and the $L^{2}_{\omega^{-1}}$ norm, and the estimate \eqref{L2errest}. The proof
of the estimate  for $u - u_{N}$ in the $H^{\alpha/2}_{\omega^{-1}}$ norm is not so straight forward.
The following lemma is helpful in establishing the $H^{\alpha/2}_{\omega^{-1}}$ error estimate.
\begin{lemma} \label{lma4Her}
Let $0 \le \mu \le 1$. For $\zeta \in H^{\mu}_{\omega}(\mrI)$, then $z \, := \, \omega \, \zeta \in 
H^{\mu}_{\omega^{-1}}(\mrI)$, with, for some $C > 0$,
\be
 \| z \|_{H^{\mu}_{\omega^{-1}}(\mrI)} \ \le \ C \,  \| \zeta \|_{H^{\mu}_{\omega}(\mrI)} \, . \label{zestH}
\ee
\end{lemma}
\textbf{Proof}:  For this proof it is convenient to use the definition of the $H^{s}_{(a , b)}(\mrI)$ spaces given by
\eqref{defHw}. \\
Let $\mu = 0$, and $\zeta \in H^{0}_{\omega}(\mrI) \, = \, L^{2}_{\omega}(\mrI)$. Then, for $z \, = \, \omega \, \zeta$
\be
\| z \|^{2}_{H^{0}_{\omega^{-1}}(\mrI)} \ = \ \| z \|^{2}_{L^{2}_{\omega^{-1}}(\mrI)} 
\ = \ \int_{\mrI} (1 - x)^{- (\alpha - \beta)} x^{- \beta} \, \big( \omega \, \zeta \big)^{2} \, dx 
\ = \  \| \zeta \|_{L^{2}_{\omega}(\mrI)} \ = \  \| \zeta \|_{H^{0}_{\omega}(\mrI)} \, .
\label{erfd0}
\ee
Next, for $\mu = 1$, let $\zeta \in C^{\infty}(\mrI) \subset H^{1}_{\omega}(\mrI)$, and let  $z \, = \, \omega \, \zeta$.
Note that $D z \, \sim \, (1 - x)^{\alpha - \beta} x^{\beta - 1} \zeta(x) \, + \, (1 - x)^{\alpha - \beta - 1} x^{\beta} \zeta(x) 
\, + \,  (1 - x)^{\alpha - \beta} x^{\beta} D\zeta(x)$, and
\begin{align}
\int_{\mrI} (1 - x)^{- (\alpha - \beta) + 1} x^{- \beta + 1}  \big( D z \big)^{2} \, dx 
&\sim \ \int_{\mrI} (1 - x)^{(\alpha - \beta) + 1} x^{\beta - 1} \zeta(x)^{2} \, dx \ + \ 
\int_{\mrI} (1 - x)^{(\alpha - \beta) - 1} x^{\beta + 1} \zeta(x)^{2} \, dx \nonumber \\
& \quad \quad   \ + \
\int_{\mrI} (1 - x)^{(\alpha - \beta) + 1} x^{\beta + 1} \, D\zeta(x)^{2} \, dx  \nonumber \\
&:= \ \mcI_{1} \ + \mcI_{2} \ + \ \| D \zeta \|^{2}_{L^{2}_{(1 - x)^{(\alpha - \beta) + 1} x^{\beta + 1}}(\mrI)} \, . \label{erfd1}
\end{align}

To bound $\mcI_{1}$ and $\mcI_{2}$ in terms of $\| \zeta \|_{L^{2}_{\omega}(\mrI)}$  and
$\| D \zeta \|_{L^{2}_{(1 - x)^{(\alpha - \beta) + 1} x^{\beta + 1}}(\mrI)}$ we use Hardy's inequality \cite[Lemma 3.2]{ber071}.
\begin{align}
\mcI_{1}
&= \int_{0}^{1/2}  (1 - x)^{(\alpha - \beta) + 1} x^{\beta - 1} \zeta(x)^{2} \, dx \ + \
 \int_{1/2}^{1}  (1 - x)^{(\alpha - \beta) + 1} x^{\beta - 1} \zeta(x)^{2} \, dx  \nonumber \\
&\lesssim \  \int_{0}^{1/2}  x^{\beta - 1} \zeta(x)^{2} \, dx \ + \
 \int_{1/2}^{1} (1 - x)^{(\alpha - \beta) + 1} x^{\beta + 1} \zeta(x)^{2} \, dx  \nonumber \\
&\lesssim \  \int_{0}^{1/2}  x^{\beta + 1} ( D \zeta(x) )^{2} \, dx \ + \ \int_{0}^{1/2}  x^{\beta + 1} \zeta(x)^{2} \, dx \ + \
 \int_{1/2}^{1} (1 - x)^{(\alpha - \beta) + 1} x^{\beta + 1} \zeta(x)^{2} \, dx  \\
 & \quad \quad \quad \quad  \ \ \mbox{(using Hardy's inequality)} \nonumber \\
&\lesssim \  \int_{0}^{1/2}  (1 - x)^{(\alpha - \beta) + 1} x^{\beta + 1} ( D \zeta(x) )^{2} \, dx \ + \ \int_{0}^{1/2}  \omega \,  \zeta(x)^{2} \, dx
 \ + \
 \int_{1/2}^{1}  \omega \,  \zeta(x)^{2} \, dx  \nonumber \\
&\lesssim \  \| \zeta \|^{2}_{L^{2}_{\omega}(\mrI)}  \ + \ \| D \zeta \|^{2}_{L^{2}_{(1 - x)^{(\alpha - \beta) + 1} x^{\beta + 1}}(\mrI)} \, .    \label{erfd2}
\end{align}
An analogous argument yields
\be
\mcI_{2} \  \lesssim \  \| \zeta \|^{2}_{L^{2}_{\omega}(\mrI)}  \ 
+ \ \| D \zeta \|^{2}_{L^{2}_{(1 - x)^{(\alpha - \beta) + 1} x^{\beta + 1}}(\mrI)} \, .   
 \label{erfd3}
\ee

Combining \eqref{erfd0}, \eqref{erfd1}, \eqref{erfd2}, and \eqref{erfd3}, we obtain
\be
  \| z \|_{H^{1}_{\omega^{-1}}(\mrI)} \ \le \ C \,  \| \zeta \|_{H^{1}_{\omega}(\mrI)} \, .
 \label{erfd4}
\ee

Estimate \eqref{erfd4} extends to $\zeta \in H^{1}_{\omega}(\mrI)$, using the density of 
$C^{\infty}(\mrI)$ in $H^{1}_{\omega}(\mrI)$.

Finally, estimate \eqref{zestH} then follows from \eqref{erfd0} and \eqref{erfd4} using interpolation. \\
\mbox{ } \hfill \qed

\begin{corollary}  \label{uerrL2}
For $f \in H^{-\alpha/2}(\mrI) \cap H^{s}_{\omega^{*}}(\mrI)$,  $s \ge -\alpha/2$, and $b$ and $c$ satisfying the 
hypothesis of Theorem \ref{thmreg13},
there
exists $C > 0$ such that for $u$ determined from
\eqref{wform1} and $u_{N}$ determined from \eqref{appx1}
\begin{align}
 \| u \, - \, u_{N} \|_{L^{2}_{\omega^{-1}}(\mrI)} \
&\le \  C \, N^{\wtilde{s} \, + \, \alpha} \, \| f \|_{H^{\wtilde{s}}_{\omega^{*}}(\mrI)} \, , \label{uL2errest}  \\
 \| u \, - \, u_{N} \|_{H^{\alpha/2}_{\omega^{-1}}(\mrI)} \ 
&\le \  C \, N^{\wtilde{s} \, + \, \alpha/2} \, \| f \|_{H^{\wtilde{s}}_{\omega^{*}}(\mrI)} \, . \label{uHerrest} 
\end{align}
\end{corollary}
\textbf{Proof}:  As commented above, \eqref{uL2errest} follows from the definition of $u_{N}$ and 
\eqref{errL2}. The estimate \eqref{uHerrest} follows from \eqref{zestH} (with $z \, = \, u - u_{N}$,
$\zeta \, = \, \phi - \phi_{N}$) and \eqref{herrest}. \\
\mbox{ } \hfill \qed

 \setcounter{equation}{0}
\setcounter{figure}{0}
\setcounter{table}{0}
\setcounter{theorem}{0}
\setcounter{lemma}{0}
\setcounter{corollary}{0}
\setcounter{definition}{0}
%
\section{Numerical Experiments}
\label{sec_num}
In this section we present three numerical experiments to
investigate the approximation of \eqref{DefProb2},\eqref{DefBC2} using
\eqref{appx1}. We compare the approximation errors with those predicted by Corollary \ref{errL2}.

For the numerical experiments we use $f(x) = 1$ and 
$f(x) = \left\{ \begin{array}{rl}
0 , &  0 < x \le 1/2 \, ,  \\
1 , &  1/2 < x < 1  
\end{array} \right. $. 
For these choices of $f$ the true solution is unknown. In order
to be able to compute a convergence rate for the approximation a very accurate approximation (using $N = 40$)
is used as the reference solution.
For the computational experiments the entries of the coefficient matrices, which require the evaluation of
integrals of weighted products of Jacobi polynomials on $\mrI$, are evaluated using the Legendre-Gauss 
quadrature rule with $200$ nodes. This ensures sufficient accuracy in order to accurately measure the
error associated with the approximation scheme \eqref{appx1}. 
We evaluate the norms of the error using the norms associated with Definition \ref{Defspace}.

The numerical convergence rate, $\kappa$, corresponding to 
$\| u_{40} \, - \, u_{N} \|_{\mbox{norm}} \lesssim N^{- \kappa}$, is presented in the tables together with the errors.
Also included are plots of the reference solution $u_{40}$, and the error $u_{40} - u_{N}$.

In Experiment 1 the data is symmetric about $x = 1/2$. However the operator is not symmetric $(r = 0.2)$, corresponding
to a preferred diffusion toward $x = 1$ over diffusion toward $x = 0$. This is reflected in the solution being
slightly skewed toward $x = 1$ (see Figure \ref{plot1}). In Experiment 2 the larger value of $r$ $(r = 0.3)$, together
with a left-to-right drift (advection) term results in a solution highly skewed to the right (see Figure \ref{plot2}).
For Experiment 3, with the diffusion and drift parameters as used in Experiment 2, the source term is taken to 
be zero for $x \in (0 , 1/2)$ and one for $x \in (1/2 , 1)$. This data results again in a solution highly skewed to
the right (see Figure \ref{plot3}).

Typically when approximating a function which is itself, or its derivative, singular at a point $x_{s}$,
the error in the approximation will be significantly larger in a neighborhood of $x_{s}$.
In the approximation scheme studied herein the correct endpoint behavior of the solution is
built into the approximation. Figures \ref{plot1}-\ref{plot3} contain plots of the error for the approximations.
In Experiments 1 and 2 the largest errors occur at the right hand endpoint, $x = 1$. Notable is that the errors 
in a neighborhood of $x = 1$ are the same order of magnitude as the errors across the interval.
For Experiment 3 the largest errors occur in a neighborhood of the discontinuity in the source term, 
around $x = 1/2$.

\textbf{Experiment 1. Fractional diffusion, reaction equation with $C^{\infty}(\mrI)$ data}. \\
For this experiment we use $\alpha = 1.60$, $r = 0.20$, $b(x) = 0$, $c(x) = 5$, and $f(x) = 1$.
Theorem \ref{thmreg11} states that even with $C^{\infty}(\mrI)$ data the regularity of the solution is bounded.
For this data $\beta = 0.93$, and 
$\wtilde{s} \, = \, \min\{ \infty, \, \alpha + (\alpha - \beta) + 1, \, \alpha +  \beta + 1 \} \, = \, 3.27$.
Corollary \ref{uerrL2} predicts that $\| u \, - \, u_{N} \|_{L^{2}_{\omega^{-1}}(\mrI)} \, \sim \, N^{-4.87}$
and $ \| u \, - \, u_{N} \|_{H^{\alpha/2}_{\omega^{-1}}(\mrI)} \, \sim \, N^{-4.07}$. The numerical
convergence rates for the errors are presented in Table \ref{table11}, and are in good agreement
with the predicted rates. A plot of the reference solution and plots of the errors are given in Figure \ref{plot1}.

\begin{table}[h!]
	\setlength{\abovecaptionskip}{0pt}
	\centering
	\caption{Experiment 1: $\alpha = 1.60$, $r = 0.20$, $b(x) = 0$, $c(x) = 5$, and $f(x) = 1$.}	\label{table11}
	\vspace{0.5em}	
	\begin{tabular}{ccccc}
		\hline
		$N$&$\|u-u_N\|_{L^2_{\omega^{-1}}}$ &$\kappa$& $\|u-u_N\|_{H^{\alpha/2}_{\omega^{-1}}}$\\
		\cline{1-5}
6&	1.05E-04&		&5.36E-04	&\\
8&	2.52E-05&	4.97 &	1.56E-04&	4.30\\ 
10&	8.62E-06&	4.81 &	6.22E-05&	4.11 \\
12&	3.61E-06&	4.77 &	2.97E-05&	4.06 \\
14&	1.74E-06&	4.76 &	1.59E-05&	4.05 \\
		\hline
		Pred.&    &4.87&     &  4.07 \\
		\hline	
	\end{tabular}
\end{table}

\begin{figure}[h!]
	\setlength{\abovecaptionskip}{0pt}
	\centering
	\includegraphics[width=3.2in,height=3.2in]{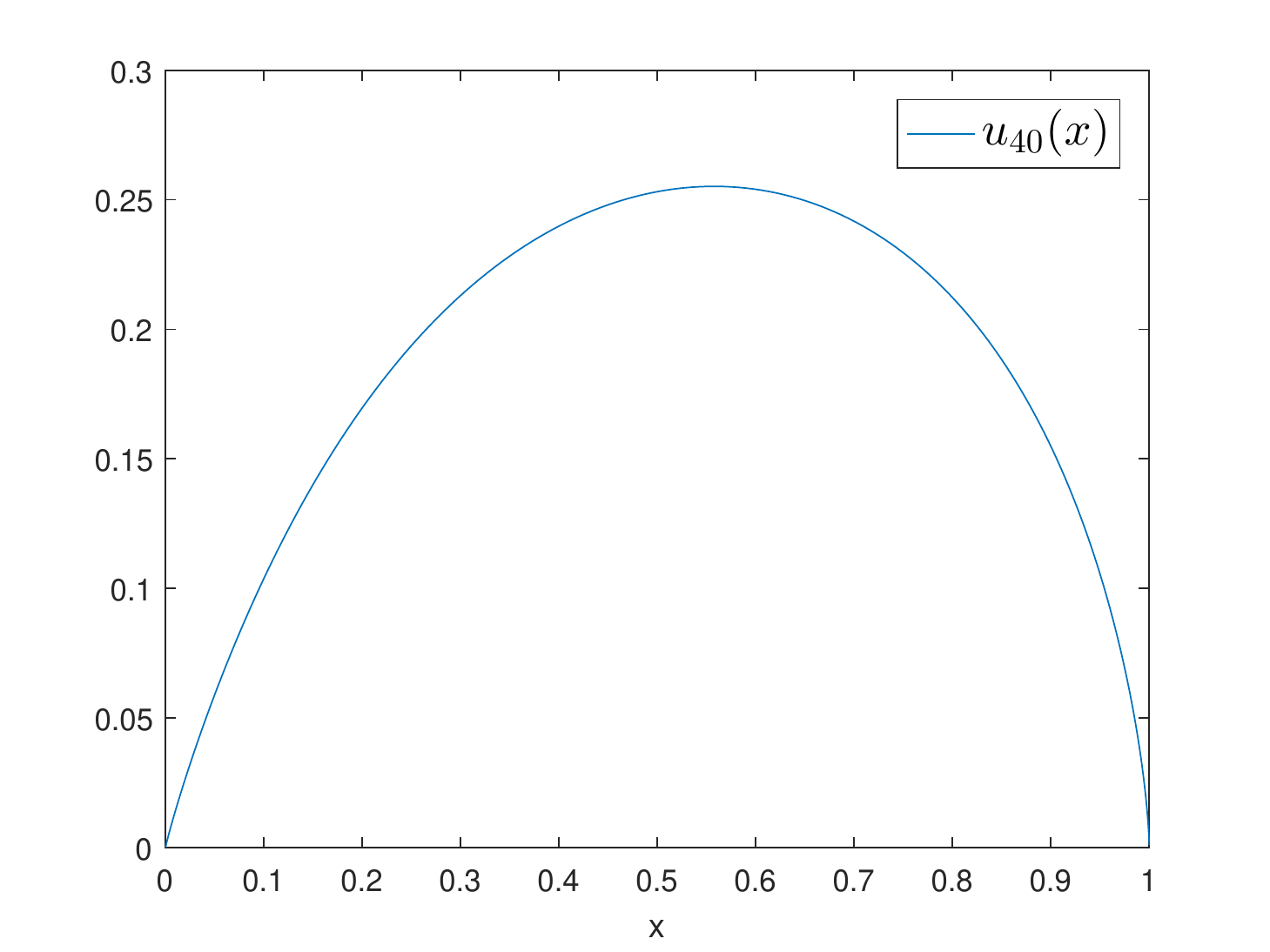}
	\includegraphics[width=3.2in,height=3.2in]{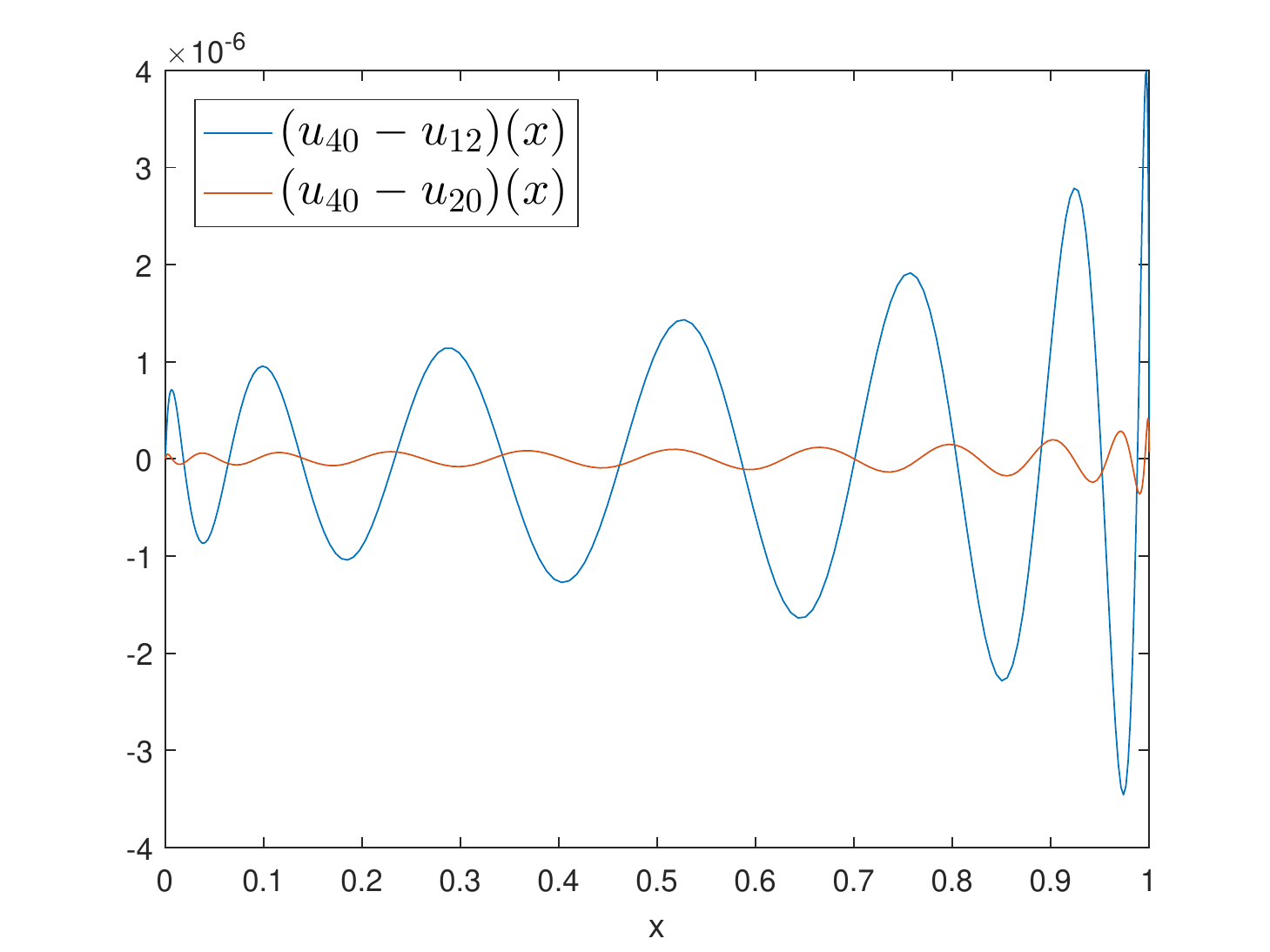}
	\caption{The plot of the reference solution $u_{40}(x)$ (left), and the plot of the errors for Experiment 1.}
	\label{plot1}
\end{figure}

\textbf{Experiment 2. Fractional diffusion, advection, reaction equation with $C^{\infty}(\mrI)$ data}. \\
For this experiment we use $\alpha = 1.40$, $r = 0.40$, $b(x) = 2$, $c(x) =  5$, and $f(x) = 1$.
As previously commented, even with $C^{\infty}(\mrI)$ data the regularity of the solution is bounded.
In addition, comparing Theorems \ref{thmreg11} and \ref{thmreg13}, the presence of an advection term
results in reduced regularity of the solution of the fractional diffusion, advection, reaction equation to that
of the fractional diffusion, advection equation. 
For this data $\beta = 0.93$, and 
$\wtilde{s} \, = \, \min\{ \infty, \, \alpha + (\alpha - \beta) - 1, \, \alpha +  \beta - 1 \} \, = \, 1.01$.
Corollary \ref{uerrL2} predicts that $\| u \, - \, u_{N} \|_{L^{2}_{\omega^{-1}}(\mrI)} \, \sim \, N^{-2.41}$
and $ \| u \, - \, u_{N} \|_{H^{\alpha/2}_{\omega^{-1}}(\mrI)} \, \sim \, N^{-1.71}$. The numerical
convergence rates for the errors are presented in Table \ref{table12}, and are in good agreement
with the predicted rates. A plot of the reference solution and plots of the errors are given in Figure \ref{plot2}

\begin{table}[h!]
	\setlength{\abovecaptionskip}{0pt}
	\centering
	\caption{Experiment 2: $\alpha = 1.40$, $r = 0.40$, $b(x) = 2$, $c(x) =  5$, and $f(x) = 1$}	\label{table12}
	\vspace{0.5em}	
	\begin{tabular}{ccccc}
		\hline
		$N$&$\|u-u_N\|_{L^2_{\omega^{-1}}}$ &$\kappa$& $\|u-u_N\|_{H^{\alpha/2}_{\omega^{-1}}}$ &$\kappa$ \\
		\cline{1-5}
12&	3.48E-03&		&2.18E-02	&\\
14&	2.47E-03&	2.21 &	1.69E-02&	1.65\\ 
16&	1.83E-03&	2.26 &	1.34E-02&	1.72 \\
18&	1.40E-03&	2.30 &	1.09E-02&	1.80 \\
20&	1.09E-03&	2.34 &	8.92E-03&	1.88 \\
		\hline	
		Pred.&    &2.41&     &  1.71 \\
		\hline	
	\end{tabular}
\end{table}

\begin{figure}[h!]
	\setlength{\abovecaptionskip}{0pt}
	\centering
	\includegraphics[width=3.2in,height=3.2in]{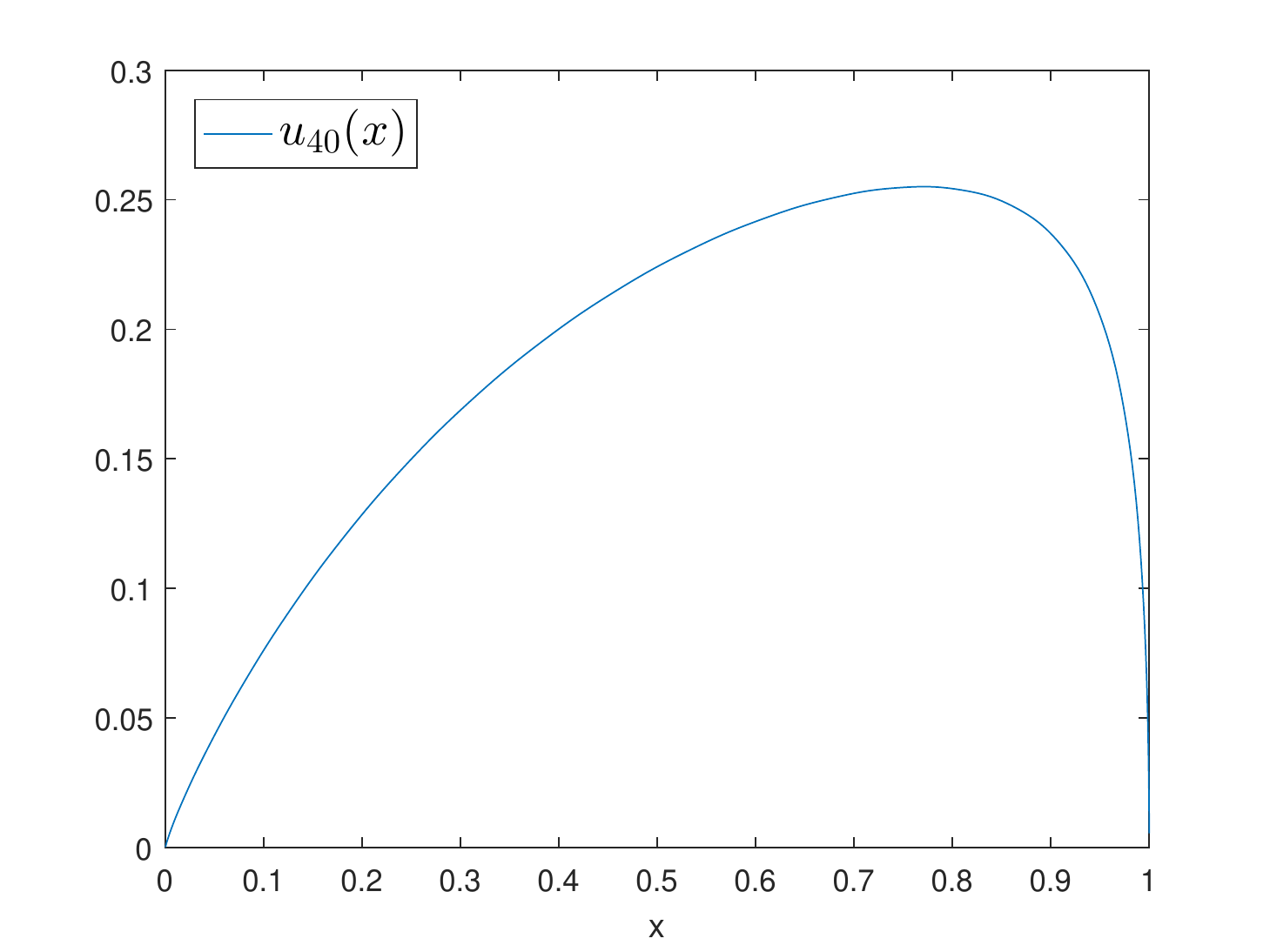}
	\includegraphics[width=3.2in,height=3.2in]{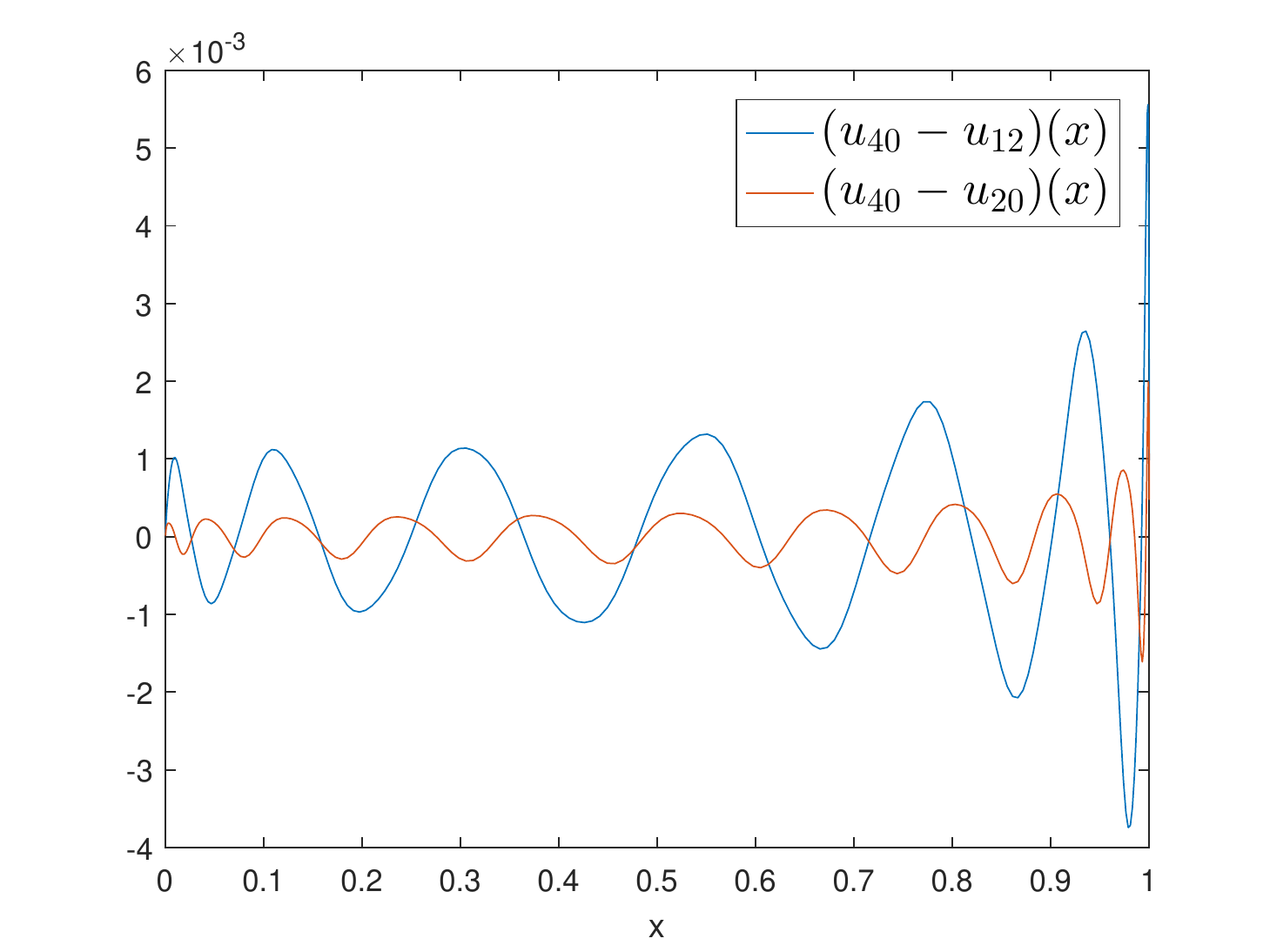}
	\caption{The plot of the reference solution $u_{40}(x)$ (left), and the plot of the errors for Experiment 2.}
	\label{plot2}
\end{figure}

\textbf{Experiment 3. Fractional diffusion, advection, reaction equation with $f \in H^{1/2 - \epsilon}_{\omega^{*}}(\mrI)$}. \\
For this experiment we use $\alpha = 1.70$, $r = 0.30$, $b(x) = 2$, $c(x) =  5$, and 
$f(x) = \left\{ \begin{array}{rl}
0 , &  0 < x \le 1/2 \, ,  \\
1 , &  1/2 < x < 1  
\end{array} \right. $.  In this case the regularity of the solution is limited by the the regularity of $f$. 
For this data $\beta = 0.91$, and 
$\wtilde{s} \, = \, \min\{ 1/2 - \epsilon, \, \alpha + (\alpha - \beta) - 1, \, \alpha +  \beta - 1 \} \, = \, 1/2 - \epsilon$.
Corollary \ref{uerrL2} predicts that $\| u \, - \, u_{N} \|_{L^{2}_{\omega^{-1}}(\mrI)} \, \sim \, N^{-2.2}$
and $ \| u \, - \, u_{N} \|_{H^{\alpha/2}_{\omega^{-1}}(\mrI)} \, \sim \, N^{-1.35}$. The numerical
convergence rates for the errors are presented in Table \ref{table333}, and are in good agreement
with the predicted rates. A plot of the reference solution and plots of the errors are given in Figure \ref{plot3}.

\begin{table}[h!]
	\setlength{\abovecaptionskip}{0pt}
	\centering
	\caption{Experiment 3: $\alpha = 1.70$, $r = 0.30$, $b(x) = 2$, $c(x) =  5$ }	\label{table333}
	\vspace{0.5em}	
	\begin{tabular}{ccccc}
		\hline
		$N$&$\|u-u_N\|_{L^2_{\omega^{-1}}}$ &$\kappa$& $\|u-u_N\|_{H^{\alpha/2}_{\omega^{-1}}}$ &$\kappa$ \\
		\cline{1-5}
	12&	3.71E-04&		&4.27E-03	&\\
	14&	2.69E-04&	2.10&	3.45E-03&	1.38\\
	16&	2.08E-04&	1.91&	2.92E-03&	1.26\\
	18&	1.61E-04&	2.18&	2.44E-03&	1.50\\
	20&	1.30E-04&	2.02&	2.11E-03&	1.40\\
		\hline	
		Pred.&    &2.20&     &  1.35 \\
		\hline	
	\end{tabular}
\end{table}

\begin{figure}[h!]
	\setlength{\abovecaptionskip}{0pt}
	\centering
	\includegraphics[width=3.2in,height=3.2in]{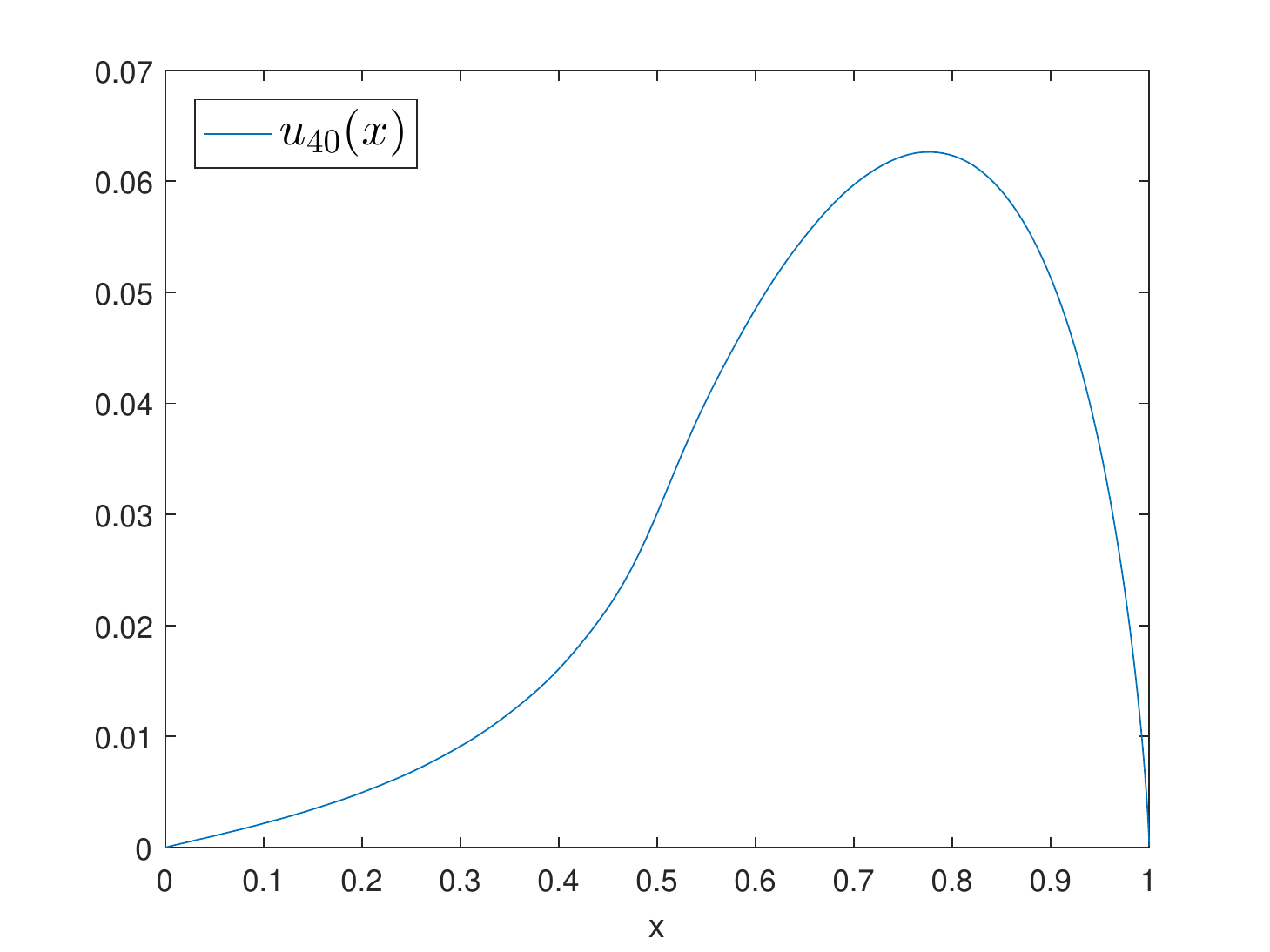}
	\includegraphics[width=3.2in,height=3.2in]{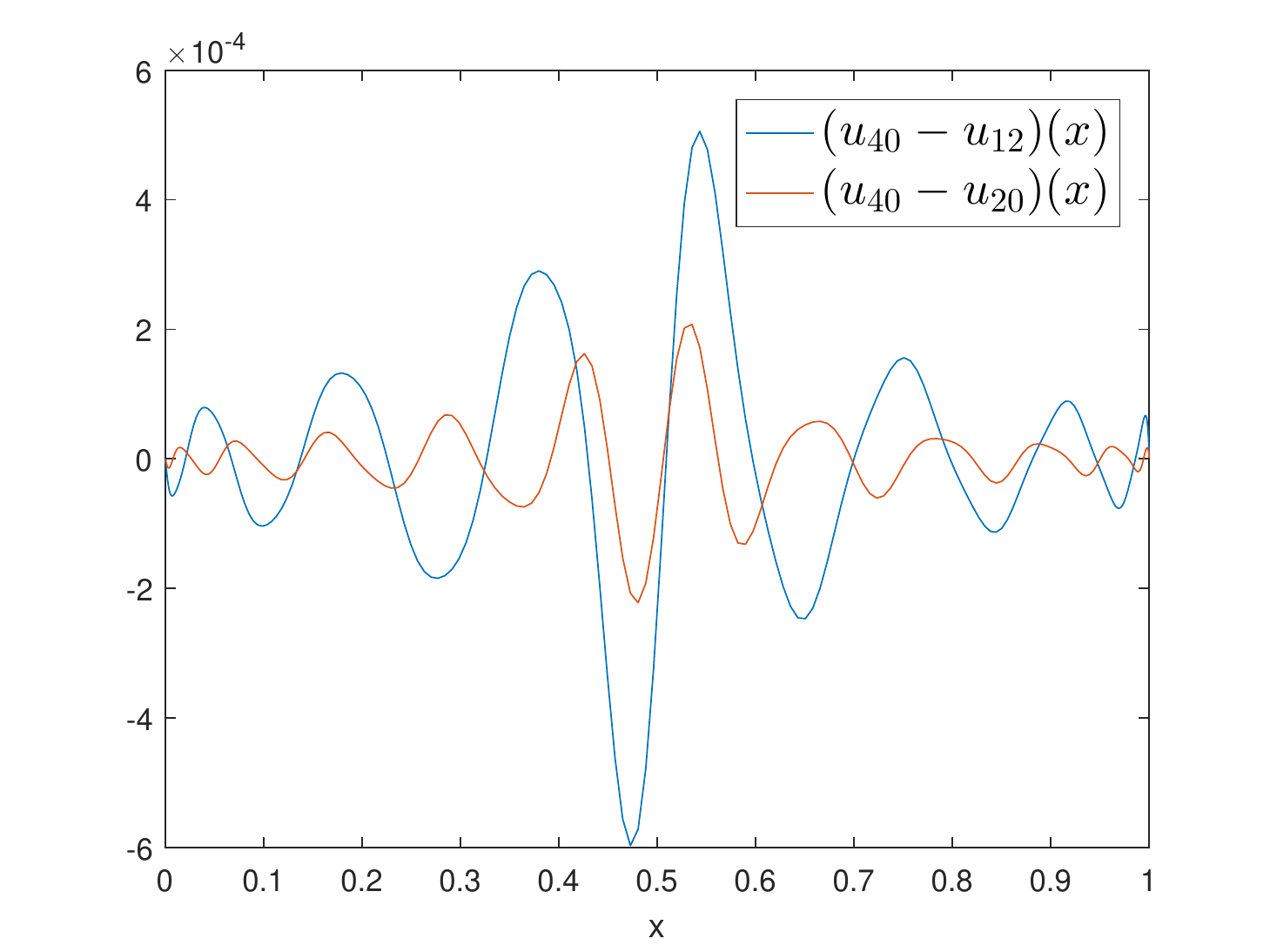}
	\caption{The plot of the reference solution $u_{40}(x)$ (left), and the plot of the errors for Experiment 3.}
	\label{plot3}
\end{figure}

\section*{Acknowledgements}
 This work was partially funded by the OSD/ARO MURI Grant W911NF-15-1-0562 and by the National Science Foundation under Grant DMS-1620194.



\end{document}